\documentclass[aop]{imsart}

\RequirePackage{amsthm,amsmath,amsfonts,amssymb,mathtools}
\RequirePackage[numbers]{natbib}
\RequirePackage[colorlinks,citecolor=blue,urlcolor=blue]{hyperref}

\usepackage[mathscr]{eucal}
\usepackage{stackengine} 
\usepackage{bm,booktabs}
\usepackage{bbm}
\usepackage[shortlabels]{enumitem}
\usepackage[usenames,dvipsnames]{xcolor}
\usepackage{comment}
\usepackage{ulem}

\usepackage{tikz}
\usepackage{pgfplots}
\usetikzlibrary{patterns}
\usepgfplotslibrary{colormaps}
\pgfplotsset{compat=1.18}
\usepackage[labelformat=simple]{subcaption}
\usepackage{caption}
\usepackage{todonotes}

\numberwithin{equation}{section} 

\newcommand{\prob}{\mathbb{P}}
\newcommand{\Prob}[1]{\prob\left(#1\right)}

\newcommand{\E}[1]{\mathbb{E}\left[ #1 \right]}
\newcommand{\En}[1]{\mathbb{E}_n\left[ #1 \right]}

\newcommand{\Var}[1]{\textup{Var}\left( #1 \right)}

\newcommand{\Cov}[1]{\textup{Cov}\left( #1 \right)}

\newcommand{\plim}{\ensuremath{\stackrel{\prob}{\longrightarrow}}}

\newcommand{\dd}{{\; \rm d}}

\newcommand{\GIRG}{\operatorname{GIRG}}
\newcommand{\GIRGn}{\GIRG^{(n)}}

\newcommand{\balpha}{\boldsymbol{\alpha}}
\newcommand{\bbeta}{\boldsymbol{\beta}}

\newcommand{\bv}{{\boldsymbol{v}}}
\newcommand{\bu}{{\boldsymbol{u}}}
\newcommand{\bw}{{\boldsymbol{w}}}
\newcommand{\bx}{{\boldsymbol{x}}}
\newcommand{\by}{{\boldsymbol{y}}}

\newcommand{\Nind}{N^{(\text{ind})}}
\newcommand{\Nsub}{N^{(\text{sub})}}
\newcommand{\Nstar}{N^{(\star)}}
\newcommand{\Jind}{J^{(\text{ind})}}
\newcommand{\Jsub}{J^{(\text{sub})}}
\newcommand{\Iind}{I^{(\text{ind})}}
\newcommand{\Isub}{I^{(\text{sub})}}
\newcommand{\Jstar}{J^{(\star)}}
\newcommand{\Istar}{I^{(\star)}}

\newcommand{\Mab}{M_{\varepsilon}^{(\balpha,\bbeta)}}
\newcommand{\M}{M_{\varepsilon}^{(\balpha^*,\bbeta^*)}}
\newcommand{\Mbar}{\overline{M}_{\varepsilon}^{(\balpha^*,\bbeta^*)}}

\newcommand{\Hk}{\mathcal{H}_k}

\DeclareCaptionSubType*{figure}

\colorlet{alpha-0}{black!50!white}
\definecolor{alpha-1}{HTML}{D7E1EE}
\definecolor{alpha-2}{HTML}{A4A2A8}
\definecolor{alpha-3}{HTML}{DF8879}
\definecolor{alpha-4}{HTML}{991F17}

\colorlet{beta-0}{black!50!white}
\colorlet{beta-1}{red}
\colorlet{beta-1}{magenta}
\colorlet{beta-2}{orange}
\definecolor{beta-3}{rgb}{0,0,1}
\definecolor{beta-4}{rgb}{0,1,0}
\definecolor{beta-5}{rgb}{1,0.5,0}
\definecolor{beta-6}{rgb}{0,1,1}
\definecolor{beta-7}{rgb}{0,0,0}

\tikzset{%
  fillcolormap/.style={/utils/exec={\pgfplotscolormapdefinemappedcolor{#1}}, draw=black, fill=mapped color}
}
\DeclareDocumentCommand\graphNode{mmm}{%
  \node[fill=alpha-#3, draw=black!70!white, very thin, circle, inner sep=0mm, minimum size=6pt] (#1) at (#2:5mm) {};
}
\DeclareDocumentCommand\graphEdge{mmm}{%
  \draw[black!70!white, ultra thick] (#1) to (#2);
  \draw[color of colormap={#3}, very thick] (#1) to (#2);
}
\DeclareDocumentCommand\graphNonEdge{mmm}{%
  \draw[color of colormap={#3}, dashed, very thick] (#1) to (#2);
}
\DeclareDocumentCommand\graphNonEdgeInduced{mmm}{%
  \draw[color of colormap={#3}, dotted, very thick] (#1) to (#2);
}

\pgfplotsset{colormap/viridis}
\tikzset{nonUnique/.style={opacity=.6,transparency group}}

\startlocaldefs
\theoremstyle{plain}

\newtheorem{theorem}{Theorem}[section]
\newtheorem{lemma}[theorem]{Lemma}
\newtheorem{proposition}[theorem]{Proposition}
\newtheorem{corollary}[theorem]{Corollary}

\endlocaldefs

\begin{document}

\begin{frontmatter}
\title{Optimal subgraphs in geometric scale-free random graphs}
\runtitle{Optimal subgraphs in geometric scale-free random graphs}

\begin{aug}
\author[A]{\fnms{Riccardo}~\snm{Michielan}\ead[label=e1]{r.michielan@utwente.nl}
\orcid{0000-0003-4642-8507}},
\author[A]{\fnms{Clara}~\snm{Stegehuis}\ead[label=e2]{c.stegehuis@utwente.nl}
\orcid{0000-0003-3951-5653}}
\author[A]{\fnms{Matthias}~\snm{Walter}\ead[label=e3]{m.walter@utwente.nl}
\orcid{0000-0002-6615-5983}}
\address[A]{University of Twente, Enschede, Netherlands\printead[presep={,\ }]{e1,e2,e3}}

\end{aug}

\begin{abstract}
Geometric scale-free random graphs are popular models for networks that exhibit as heavy-tailed degree distributions, small-worldness and high clustering. In these models, vertices have weights that cause the heavy-tailed degrees and are embedded in a metric space so that close-by groups of vertices tend to cluster. The interplay between the vertex weights and positions heavily affects the local structure of the random graph, in particular the occurrence of subgraph patterns, but the dependencies in these structures and weights make them difficult to analyze. In this paper we investigate subgraph counts using a \textit{divide et impera} strategy: first counting the number of subgraphs in specific classes of vertices; then computing which class yields maximum contribution. Interestingly, the scaling behavior of induced and general subgraphs in such geometric heavy-tailed random graphs is closely related to the solution of a mixed-integer linear program which also shows that subgraphs appear predominantly on vertices with some prescribed degrees and inter-distances. Finally, we derive precise asymptotics for trees and Hamiltonian subgraphs.
\end{abstract}



\end{frontmatter}


\section{Introduction}

Network subgraphs have been an object of intense study over the past decades, both from a combinatorial perspective, but more recently also in analyzing real-world network data. Here, subgraph counts have gained attention through so-called motifs: statistically significant subgraphs of networks~\cite{milo2002}. The typical setting used to analyze motifs is to compare the frequency of a particular chosen subgraph in a network with that in an associated network null-model or random graph model. This investigation initially started in social networks \cite{holland1976,holland1977}, biological networks \cite{shen2002,ptacek2005,ma2005,alon2007} and ecological networks \cite{connor1979}, but the interest in network motifs rapidly expanded to other fields including finance and communication science \cite{saracco2015,saracco2016,jiang2013}. 

One of the main criticisms of the motif analyses is that the null-model compared to which statistical significance is claimed is often unrealistic. Indeed, in most examples, a networked data set is compared to the configuration model~\cite{molloy1995}: a random graph model that only matches the exact degree sequence of the network. In such random graphs, any type of non-tree subgraph is rare, in contrast to most real-world networks which are often heavily clustered. Therefore, comparing to the configuration model often yields that every subgraph containing a cycle appears more frequently in the data set than in the random graph model, while such structures are not a rare property but rather a commonality of real-world data.
Investigating small subgraphs in more realistic random graph models is therefore of vital importance to better understand which specific subgraph frequencies emerge naturally from realistic network features.

Recurrent properties in large real-world networks are heavy-tailed degree sequences, small distances and large clustering coefficient. Several recent random models create graphs with all these properties by using vertex weights to tune the degree sequence to a power-law distribution, and by introducing an underlying geometry that can explain the clusteredness. Examples of such models are hyperbolic random graphs \cite{krioukov2010}, spatial preferential attachment models \cite{jacob2015}, and spatial inhomogeneous random graphs \cite{bringmann2019,deprez2019,gracar2022}. The local limits of this class of random graphs deviates substantially from the usual locally tree-like structure of many complex network models \cite{vdhofstad2016,vdhofstad2023}. Therefore, the occurrence of small subgraphs differs significantly when compared to random graphs where geometric features are absent. However, the geometric setting of these random graph models make them mathematically difficult to analyze, as the geometry ruins edge independence. Furthermore, the power-law degrees create so-called degree-degree correlations, which destroy the product form of the edge probabilities.

For this reason, subgraph counts have mainly been studied in random graphs with more homogeneous degree distribution and no underlying geometric structure \cite{ruciniski1985,nowicki1988,itzkovitz2003}. More recently, a novel powerful technique has been introduced to deal with degree inhomogeneity in random graphs for counting subgraph occurrences \cite{stegehuis2019,vdhofstad2021}. This approach allows one to obtain a precise estimate of any subgraph count and at the same time to locate where within the network the subgraph appears most frequently. However, this technique strongly relies on independent edge presences, which do not hold for geometric random graph models. 

In this paper, we investigate the number of subgraphs in a model that combines scale-freeness and geometry, the geometric inhomogeneous random graph \cite{bringmann2019}. We show that most of the subgraphs emerge from set of vertices with specific degrees and inter-distances. We create a new technique to overcome the difficulties in analyzing dependent, geometric edge presences by introducing and analyzing a mixed-integer linear program (MILP) that is able to capture the edge dependencies in its constraints. We use its optimizer to prove the scaling of all subgraph counts.  The MILP not only determines the scaling of the subgraph counts, it also gives insights into the typical locations of all possible subgraphs, in terms of the vertex degrees and their distances.
We show that via a sequence of related MILPs, we can also deduce the uniqueness and possible symmetries in these optimal locations.
To our knowledge, this is one of the first papers that uses mixed integer linear programs to analyze random graphs, demonstrating that they can be a powerful tool to show the dominant behavior of complex random graph models. Furthermore, the method is general, and applies to all possible fixed subgraphs, while to the best of our knowledge, all previous results are restricted to specific subgraphs with extreme symmetry properties, such as cliques~\cite{michielan2021,Oh2023,Blasius2017,Blaesius2023, Soenmez2024}. 

Our methods show that in these complex random graph models, subgraphs often appear at very specific degrees and distances. This highlights the intricate interplay between power-law degrees and geometry: some subgraphs are typically formed by close by vertices of low degrees caused by geometry, while other subgraphs are dominated by high-degree vertices that are typically not close by in the geometry.

\subsection{Notation}
Throughout this paper we build a random graph model on the set of vertices $[n]= \{1,2,\dots, n\}$ and study the large $n$ limit. We therefore use classical asymptotic notation, in terms of the graph size $n$. For any two non-negative functions $f(n),g(n)$ we will write: $f(n)=o(g(n))$ if $\lim_{n \to \infty} f(n)/g(n) = 0$; $f(n) = O(g(n))$ if $\limsup_{n \to \infty} f(n)/g(n) < \infty$; $f(n) = \Omega(g(n))$ if $\liminf_{n \to \infty} f(n)/g(n) > 0$; $f(n) = \Theta(g(n))$ if $f(n) = O(g(n))$ and $f(n) = \Omega(g(n))$. Moreover, we will say that a sequence of events $\{\mathcal{E}_n\}_{n \geq 1}$ happens with high probability (w.h.p.) if $\lim_{n \to \infty} \Prob{\mathcal{E}_n}=1$. Finally, we will denote by $X_n \plim X$ the convergence in probability of a sequence of random variables, meaning that $\lim_{n\to \infty} \Prob{|X_n - X|>\varepsilon} = 0$, for all $\varepsilon>0$.

\subsection{Model}

We now describe the geometric inhomogeneous random graphs (GIRGs) \cite{bringmann2019}. In this class of models, the vertex set is $V = [n]:=\{1,...,n\}$ and each vertex $i$ is associated with a weight $w_i>0$. The weight sequence follows a power law distribution with exponent $\tau > 0$. We will assume that $\bw = \{w_i\}_{i \in [n]}$ are i.i.d. samples of a same random variable $W$ such that
\begin{equation}\label{eq:weightdistribution}
    \Prob{W > x} = c x^{1-\tau} (1 + o(1)), 
\end{equation} 
for all $x \geq w_0$, for some global constants $w_0,c>0$. In this paper, we assume that $\tau \in (2,3)$, that is, the weights have finite mean $\mu$, but infinite variance. For simplicity, we assume that $w_0 = 1$. Each vertex $i$ is also assigned a position $x_i$. As ground space, we consider the $d$-dimensional torus $\mathbb{T}^d = \mathbb{R}^d/\mathbb{Z}^d$, which can be described as the $d$-dimensional cube $[0,1]^d$ where opposite boundaries are identified, endowed with the infinity norm distance
\begin{equation*}
    ||x-y|| := \max_{j=1,...,d} \min \left\{ |x^{(j)} - y^{(j)}|, 1 - |x^{(j)} - y^{(j)}| \right\},
\end{equation*}
for any $x,y \in \mathbb{T}^d$. The vertex positions $\bx = \{x_i\}_{i \in [n]}$ are sampled i.i.d uniformly over $\mathbb{T}^d$.

Any two vertices $u,v \in [n]$ are then independently connected with probability $p_{uv}$ determined by the weights and positions of the vertices
\begin{equation} \label{eq:edgeprob}
    p_{uv} = p(w_u,w_v,x_u,x_v) = 1 \wedge \left(\frac{w_u w_v}{n\mu||x_u-x_v||^{d}}\right)^\gamma,
\end{equation}
where $\mu = \mathbb{E}[W]$ and $\gamma > 1$ is a fixed parameter. The choice $\gamma = \infty$ is also allowed, and in this case
\begin{equation} \label{eq:edgeprobtresh}
    p_{uv} = p(w_u,w_v,x_u,x_v) = \mathbbm{1}_{\left\{w_u w_v > n\mu||x_u-x_v||^{d}\right\}}.
\end{equation}
We denote the resulting random graph by $\GIRG^{(n)}$, where dependency on the parameters $\tau,\gamma,d$ is implicit. The connection probability in \eqref{eq:edgeprob}--\eqref{eq:edgeprobtresh} reflects the idea that nearby (similar) vertices are likely to be connected, generating a geometric structure. At the same time, high-degree (popular) vertices are also likely to connect. 

Sometimes, it will be useful to consider a \textit{cutoff} version of the $\GIRG$, that is, conditioning $\GIRGn$ on the realization of the event
\begin{equation}\label{eq:cutoff}
    \mathcal{E}_{\bar{\varepsilon}} := \left\{\max_{v \in V} w_v < n^{1/(\tau-1)}/\bar{\varepsilon}, \min_{u,v \in V} ||x_u-x_v|| > \bar{\varepsilon} n^{-1/d} \right\},
\end{equation}
for some $\bar{\varepsilon}>0$. It is immediate to see that for any $\delta > 0$, there exists $\bar{\varepsilon}$ such that $\mathbb{P}(\mathcal{E}_{\bar{\varepsilon}}) > 1-\delta$, meaning that $n^{1/(\tau-1)}$ and $n^{-1/d}$ are the natural cutoff for vertex weights and inter-distances in $\GIRGn$.

\subsection{Subgraph counting}
In this paper we analyze the number of copies of fixed subgraphs in $\GIRGn$ when the number of vertices $|V|=n$ tends to infinity, while the size $k$ of the subgraph is fixed. 
As $\GIRGn$ is a random graph, the number of occurrences of a specified subgraph is a random variable, and we will show that it obeys a specific scaling in $n$, with high probability. We now describe the problem more precisely, while introducing some notation.

Let $H=(\mathcal{V},\mathcal{E})$ be a deterministic graph, with $\mathcal{V}=[k]$ and $k=O(1)$. We denote by $\Nind(H)$ the number of induced subgraphs of $\GIRGn$ that are equal to $H$, and by $\Nsub(H)$ the number of occurrences of $H$ as general subgraph in $\GIRGn$. That is,
\begin{align}
    \Nind(H) = \sum_{\bv \in V_k} \mathbbm{1}_{\{H = \GIRGn|_{\bv}\}}, \label{eq:Nind} \\
    \Nsub(H) = \sum_{\bv \in V_k} \mathbbm{1}_{\{H \subseteq \GIRGn|_{\bv}\}}, \label{eq:Nsub}
\end{align}
where $V_k = \{(v_1,...,v_k) \subset V\}$ is the collection of all ordered sets of vertices of size $k$ in $V$ and $\GIRGn|_{\bv}$ denotes the induced subgraph of $\GIRGn$ on the ordered vertex set $\bv$. The relation $H = \GIRGn|_{\bv}$ means $(v_i,v_j) \in E$ if and only if $(i,j) \in \mathcal{E}$, whereas $H \subseteq \GIRGn|_{\bv}$ means $(v_i,v_j) \in E$ whenever $(i,j) \in \mathcal{E}$. 

The expected number of copies of $H$ in $\GIRG^{(n)}$ can be written as
\begin{equation}\label{eq:expNind}
\begin{split}
    \E{\Nind(H)} &= \sum_{\bv \in V_k} \Prob{\GIRGn|_{\bv}=H} 
    \\
    &= \frac{n!}{(n-k)!} \int \Prob{ \GIRGn(\bw,\bx)|_{\bv}=H} \dd \mathbb{P}(\bw,\bx).
\end{split}
\end{equation}
A similar expression for holds for $\Nsub(H)$, replacing $H = \GIRGn|_{\bv}$ with $H \subseteq \GIRGn|_{\bv}$. When $H$ is a $k$-clique, the integral in \eqref{eq:expNind} can be computed explicitly and the expected number of $k$-cliques is asymptotically linear in $n$ or scales as $n^{k(3-\tau)/2}$, depending on $k$ and $\tau$ \cite{michielan2021}. However, for general $H$ there is no existing method to solve such integrals because the weights and spatial coordinates make the edge presences of $H$ dependent. 

Our main step to overcome this problem is to first focus on counting the number of subgraphs that appear on specific subsets of $V_k$ only. Namely, if $A \subseteq V_k$, we define
\begin{align}\label{eq:Nha}
    \Nind(H;A) = \sum_{\bv \in A} \mathbbm{1}_{\{H = \GIRGn|_{\bv}\}}, \quad     \Nsub(H;A) = \sum_{\bv \in A} \mathbbm{1}_{\{H \subseteq \GIRGn|_{\bv}\}},
\end{align}
In the next section, we will show that first computing $\Nind(H;A)$ and $\Nsub(H;A)$ generates a general approach to estimate all fixed-size subgraph counts $H$, addressing the subgraph counting problem systematically.

\subsection{Organization of the paper}
We present our main results in Section \ref{sec:results}, where we state asymptotic results for subgraph counts and introduce an optimization problem, together with its mixed-integer linear program reformulation.  
We then provide the proofs. In Section \ref{sec:proofgenericsubgraphs} we prove the existence of \textit{optimal} regions in the GIRG, where the subgraph count is maximal. Finally, in Section \ref{sec:proofparticularsubgraphs} we prove sharp results for specific classes of subgraphs.

\section{Results}\label{sec:results}

The integral in \eqref{eq:expNind} is difficult to solve in general, as the edge presences of $H$ are dependent. Thus, we apply a \textit{divide et impera} strategy: first, we count subgraphs formed on vertices with prescribed weights and inter-positions; then, we derive and solve an optimization problem to determine where the subgraphs are more likely to appear.

More precisely, for any $\zeta \in \mathbb{R}$ and $\varepsilon \in (0,1)$ we define the interval $I_{\varepsilon}(n^{\zeta}) = [\varepsilon n^{\zeta}, n^{\zeta}/\varepsilon]$. For any $\balpha= \{\alpha_i\}_{i \in [k]}$ and $\bbeta= \{\beta_{\{i,j\}}\}_{i,j \in [k], i \neq j}$, we define the set $\Mab \subseteq V_k$ as 
\begin{equation}\label{eq:Mab}
\Mab = \left\{(v_1,...,v_k) \subset V: w_{v_i} \in I_{\varepsilon}(n^{\alpha_i}) \; \forall i, ||x_{v_i} - x_{v_j}|| \in  I_{\varepsilon}(n^{\beta_{\{i,j\}}}) \; \forall i,j \right\}.
\end{equation}
The random set $M_{\varepsilon}^{(\balpha,\bbeta)}$ contains all the lists of $k$ vertices in $\GIRG^{(n)}$ such that their weights scale as $n^{\alpha_1}\dots,n^{\alpha_k}$, and their distances on $\mathbb{T}^d$ scale as $n^{\beta_{\{i,j\}}}$ for any $i \neq j$. \\
Since weights in $\GIRGn$ are independently sampled from \eqref{eq:weightdistribution}, for any $\delta>0$ there exists $C>0$ such that
\begin{equation*}
    \limsup_{n \to \infty} \Prob{\max_{v \in V} w_v > C n^{1/(\tau-1)}} < \delta.
\end{equation*}
Thus, whenever $\alpha_i > \frac{1}{\tau-1}$ for any $i \in [k]$, $\Mab$ is empty with high probability.
Similarly, as the random positions are independently and uniformly sampled from $\mathbb{T}^d$, for any $\delta>0$ there exists $c>0$ such that
\begin{equation*}
    \limsup_{n \to \infty} \Prob{\min_{u,v \in V} ||x_u - x_v|| < c n^{-1/d}} < \delta,
\end{equation*}
so that $\Mab$ is empty with high probability whenever $\beta_{\{i,j\}} < -\frac{1}{d}$ for some $i,j$.
Moreover, the $\bbeta$ variables need to satisfy structural constraints that emerge from the torus norm. Indeed, if $||x_i-x_j|| \propto n^{\beta_{\{i,j\}}}$ and $||x_j-x_k|| \propto n^{\beta_{\{j,k\}}}$, then $||x_i-x_k|| = O(n^{\max(\beta_{\{i,j\}},\beta_{\{j,k\}})})$, as otherwise the triangle inequality would imply that the distance from $i$ to $j$ to $k$ would be shorter than the distance between $i$ and $k$. In other words, $\bbeta$ should satisfy satisfy $\beta_{\{i,j\}}\leq \max(\beta_{\{i,s\}},\beta_{\{j,s\}}), \forall s\neq i,j$. Thus, we define the \textit{feasible region}
\begin{multline}\label{eq:feasibleregion}
    \mathcal{F}:=\Big\{ \left(\{\alpha_i\}_{i\in[k]}, \{\beta_{\{i,j\}}\}_{i,j \in [k]} \right) : \alpha_i \in [0,1/(\tau-1)], \beta_{\{i,j\}} \in [0,-1/d], \\
    \beta_{\{i,j\}}\leq \max(\beta_{\{i,s\}},\beta_{\{j,s\}}),\forall i,j,s \in [k] \Big\},
\end{multline}
so that $\Mab$ is empty with high probability whenever $(\balpha,\bbeta) \not \in \mathcal{F}$.
In particular, observe that conditionally on the event $\mathcal{E}_{\bar{\varepsilon}}$ defined in \eqref{eq:cutoff}, 
\begin{equation*}
    \Nsub(H) = \Nsub\Big(H; \bigcup_{(\balpha,\bbeta) \in \mathcal{F}} M_{\bar{\varepsilon}}^{(\balpha,\bbeta)}\Big) \quad \text{and} \quad \Nind(H) = \Nind\Big(H; \bigcup_{(\balpha,\bbeta) \in \mathcal{F}} M_{\bar{\varepsilon}}^{(\balpha,\bbeta)}\Big)
\end{equation*}
almost surely.

\subsection{Optimal subgraphs}\label{sec:optimalsubgraphs}
We now present our results on subgraphs with specific vertex parameters. 
We denote the occurrences of $H$ as induced or general subgraph in $\GIRGn$ with vertices in $M_{\varepsilon}^{(\balpha,\bbeta)}$ by $\Nind(H;\Mab)$ and $\Nsub(H;\Mab)$ as in Equation~\eqref{eq:Nha}.
Our first lemma shows how the expected subgraph counts in $\Mab$ scales in $n$:

\begin{lemma}\label{lemma:expNHMab}
For any $H=(\mathcal{V},\mathcal{E})$ of size $k$ and $(\balpha,\bbeta) \in \mathcal{F}$, define
\begin{multline}
\label{eq:fh}
    f_H(\balpha,\bbeta) := k + (1 - \tau) \sum_{i = 1}^{k} \alpha_i + d \sum_{j=2}^{k} \min_{i:i<j} \beta_{\{i,j\}} \\ + \gamma \sum_{\substack{i,j \in [k] \\ (i,j) \in \mathcal{E}}} \min(\alpha_i + \alpha_j - d \beta_{\{i,j\}} - 1, 0).
\end{multline}
Then, for any fixed $\varepsilon \in (0,1)$,
\begin{equation*}
    \E{\Nsub(H;\Mab)} = \Theta\left(n^{f_H(\balpha,\bbeta)}\right),
\end{equation*}
and
\begin{equation*}
    \E{\Nind(H;\Mab)} = 
    \begin{cases}
    \Theta\left(n^{f_H(\balpha,\bbeta)}\right) & \text{if $\alpha_i + \alpha_j \leq d \beta_{\{i,j\}} + 1, \forall (i,j) \not \in \mathcal{E}$,}\\
    0 & \text{otherwise.}
    \end{cases}
\end{equation*}
\end{lemma}

The proof of this result can be found in Section \ref{sec:proofgenericsubgraphs}. Lemma \ref{lemma:expNHMab} shows that, for different choices of $(\balpha,\bbeta)$, the sets $\Mab$ yield different contributions to the expected number of copies of $H$ in $\GIRGn$. Next, the natural question is which $\Mab$ provides the maximum contribution. To this end, we formulate an optimization problem for the expected number of general subgraphs and find $(\balpha,\bbeta)$ that maximize $f_H(\balpha,\bbeta)$:
\begin{subequations}
  \label{eq:optproblem}
  \begin{alignat}{7}
    & \text{max }
      & k + (1 &- \tau) \sum_{i=1}^{k} \alpha_i + d \sum_{j=2}^{k} \min_{i:i<j} \beta_{\{i,j\}} &\;+\;& \gamma \hspace{-1ex} \sum_{\substack{i,j \in [k] \\ (i,j) \in \mathcal{E}}} \hspace{-1ex} \min(\alpha_i + \alpha_j - d \beta_{\{i,j\}} - 1, 0) \quad \label{eq:optproblem_objective} \\
    & \text{over }
      & \alpha_i    &\in [0,1/(\tau-1)] \label{eq:optproblem_alpha_boundary}
        &\quad& \forall i \in [k] \\
    &  & \beta_{\{i,j\}}  &\in [-1/d,0]  \label{eq:optproblem_beta_boundary}
        &\quad& \forall i,j \in [k] : i < j \\
    & \text{s.t. }
      & \beta_{\{i,j\}} &\leq \max(\beta_{\{i,s\}},\beta_{\{s,j\}})
        &\quad& \forall i,j,s \in [k] : i < j,~ s \neq i,~ s \neq j.\label{eq:optproblem_normal_constraint}
  \end{alignat}
\end{subequations}
For the expected number of induced subgraphs, we formulate the same optimization problem as in \eqref{eq:optproblem}, including the extra set of constraints:
\begin{equation}\label{eq:optproblem_extra_constraint}
    \alpha_i+\alpha_j \leq 1 + d \beta_{\{i,j\}}\quad \forall (i,j) \not \in \mathcal{E}
\end{equation}
It is easily checked that the optimization problems~\eqref{eq:optproblem} and ~\eqref{eq:optproblem}+\eqref{eq:optproblem_extra_constraint} are essentially independent of the dimension parameter $d$:

\begin{proposition}
  \label{thm:optproblem_dimension_independent}
  A vector $(\balpha,\bbeta)$ is feasible or optimal for~\eqref{eq:optproblem} or satisfies~\eqref{eq:optproblem_extra_constraint} for a particular dimension $d = \bar{d}$ if and only if the vector $(\balpha, \bar{d} \cdot \bbeta)$ satisfies the respective property for $d = 1$.
\end{proposition}
\begin{proof}
    Let $d\geq1$. Under the transformation $\tilde{\bbeta} = d \cdot \bbeta$, $f_{H}^{(d)}(\balpha,\bbeta) = f_{H}^{(1)}(\balpha,\tilde{\bbeta})$, where $f_{H}^{(d)}$ denotes the objective function~\eqref{eq:fh} for a given $d$. Furthermore, we can rewrite~\eqref{eq:optproblem_beta_boundary}--\eqref{eq:optproblem_extra_constraint} as
    \begin{align*}
        \tilde{\beta}_{\{i,j\}} \in [-1,0], &\qquad \forall i,j\in[k]: i<j, \\
        \alpha_i + \alpha_j \leq 1 + \tilde{\beta}_{\{i,j\}}, &\qquad \forall (i,j) \not \in \mathcal{E},
    \end{align*}
which are the same constraints as for the optimization problem with $d=1$.
\end{proof}

Another interesting property is that $\gamma$ does not play any role in the optimization problem for general subgraphs~\eqref{eq:optproblem}. Indeed, general subgraphs are always formed to minimize the \textit{edge-energy}:
\begin{proposition}
  \label{thm:optproblem_gamma_independent}
  The vector $(\balpha^*,\bbeta^*)$ is optimal for~\eqref{eq:optproblem} at a given $\gamma>1$ if and only if it is optimal for $\gamma = \infty$, that is, if and only if $(\balpha^*,\bbeta^*)$ maximizes the function $$f(\balpha,\bbeta):= k + (1 - \tau) \sum_{i = 1}^{k} \alpha_i + d \sum_{j=2}^{k} \min_{i:i<j} \beta_{\{i,j\}}$$
  subject to \eqref{eq:optproblem_alpha_boundary}, \eqref{eq:optproblem_beta_boundary}, \eqref{eq:optproblem_normal_constraint} and
  \begin{equation}
      \alpha_i + \alpha_j \geq 1 + d \beta_{\{i,j\}}, \qquad \forall (i,j) \in \mathcal{E}.
  \end{equation}
\end{proposition}

The proof of Proposition~\ref{thm:optproblem_gamma_independent} can be found in Section \ref{sec:proofgenericsubgraphs}. As an immediate consequence of Lemma \ref{lemma:expNHMab} we get that the expected number of subgraphs is dominated by contributions from the optimal solution $(\balpha^*,\bbeta^*)$:

\begin{corollary}\label{thm:subgraphGIRG}
Let $H$ be a subgraph of size $k$. Let $\star \in \{\text{ind},\text{sub}\}$ denote the induced or general subgraph setting. Suppose that the corresponding problem \eqref{eq:optproblem} or \eqref{eq:optproblem}$+$\eqref{eq:optproblem_extra_constraint} has unique solution $(\balpha^*,\bbeta^*)$. Then, for any $(\balpha,\bbeta)\in\mathcal{F} \setminus \{(\balpha^*,\bbeta^*)\}$ and $\varepsilon \in (0,1)$,
\begin{equation*}
    \lim_{n \to \infty}\frac{\E{N^{(\star)}\left(H;M_{\varepsilon}^{(\balpha,\bbeta)}\right)}}{\E{N^{(\star)}\left(H;M_{\varepsilon}^{(\balpha^*,\bbeta^*)}\right)}} = 0.
\end{equation*}
\end{corollary}

\begin{proof}
Since $(\balpha^*,\bbeta^*)$ is the unique solution of \eqref{eq:optproblem} (or \eqref{eq:optproblem}+\eqref{eq:optproblem_extra_constraint}) it follows that 
\begin{equation*}
    f_H(\balpha,\bbeta) < f_H(\balpha^*,\bbeta^*),
\end{equation*}
for any $(\balpha,\bbeta) \in \mathcal{F}$ with $(\balpha,\bbeta) \neq (\balpha^*,\bbeta^*)$.
Then, by Lemma \ref{lemma:expNHMab},
\begin{equation*}
    \frac{\E{N^{(\text{sub})}(H;\Mab)}}{\E{N^{(\text{sub})}(H;\M)}} = \frac{\Theta(n^{f_H(\balpha,\bbeta)})}{\Theta(n^{f_H(\balpha^*,\bbeta^*)})} \longrightarrow 0, \quad \text{as $n \to \infty$.}
\end{equation*}
The case of induced subgraphs follows similarly.
\end{proof}

Moreover, we obtain a stronger concentration result for subgraphs $H$ that achieve their solution on particular values $\balpha$ and $\bbeta$.

\begin{theorem}\label{thm:subgraphGIRGprob}
Under the assumptions of Corollary \ref{thm:subgraphGIRG}, suppose $(\balpha^*,\bbeta^*) \equiv (1/2, 0)$ or $(\balpha^*,\bbeta^*) \equiv (0, -1/d)$. Then, for any fixed $\varepsilon \in (0,1)$,
\begin{equation}\label{eq:subgraphGIRGprobbound}
    \lim_{n \to \infty}\Prob{
        g_1(\varepsilon) \leq
        \frac{N^{(\star)}(H;\M)}{n^{f_H(\balpha^*,\bbeta^*)}} 
        \leq g_2(\varepsilon)
        } = 1 
\end{equation}
for some positive functions $g_1,g_2$ not depending on $n$. Moreover, if $(\balpha^*,\bbeta^*)$ is the unique solution of the optimization problem, then for any $(\balpha,\bbeta)\in\mathcal{F} \setminus \{(\balpha^*,\bbeta^*)\}$ and $\varepsilon \in (0,1)$,
\begin{equation*}
    \frac{N^{(\star)}(H;\Mab)}{N^{(\star)}(H;\M)}\plim 0.
\end{equation*}
\end{theorem}

The proof of this result can be found in Section \ref{sec:proofgenericsubgraphs}. Observe that, under the hypothesis of Theorem \ref{thm:subgraphGIRGprob}, these particular subgraph counts are dominated by the contribution over set of vertices either with degree $\propto \sqrt{n}$ or at distance $\propto n^{1/d}$ from each other. Therefore, in the first case the leading contribution emerges regardless of geometry and is caused by the power-law weights. In the second case, the spatial embedding of the vertices plays a fundamental role, but the power-law weights do not. For this reason, we identify two distinct phases that we call \textit{geometric} and \textit{non-geometric}. 

Figures \ref{fig_k4_taulow_subgraph}, \ref{fig_k4_tauhigh_subgraph}, \ref{fig_k5_taulow_subgraph} and \ref{fig_k5_tauhigh_subgraph} show that the solution to the optimization problem for general subgraphs is uniquely achieved on $(\balpha,\bbeta) \equiv (1/2,0)$ or $(0,-1/d)$ for a large class of subgraphs (all subgraphs with only light gray nodes and dark blue edges, or only orange nodes and yellow edges). This class contains mostly relatively dense subgraphs, and contains for example all cliques. The same happens for induced subgraphs, as shown in Figures \ref{fig_k4_taulow_induced}, \ref{fig_k4_tauhigh_induced}, \ref{fig_k5_taulow_induced} and \ref{fig_k5_tauhigh_induced}. In particular, these dense subgraphs manifest a geometric phase when $\tau$ is large, whereas a non-geometric phase for a smaller value $\tau$.

\subsection{Total subgraph counts}\label{sec:totalsubgraphs}

The results in the previous section focus on the occurrence of $H$ on sets of the kind $\Mab$, given $(\balpha,\bbeta) \in \mathcal{F}$. We now deal with the total number of occurrences $\Nsub(H)$ and $\Nind(H)$. 
First, an immediate application of Markov inequality shows that the number of occurrences of any $H$ in $\GIRGn$ grows at least linearly or super-linearly with $n$ (depending on $k$ and $\tau$), with high probability:

\begin{corollary}
\label{thm:lineargrowth}
Let $H$ be a subgraph of size $k$ and $\star \in \{\text{ind},\text{sub}\}$.
\begin{itemize}
    \item When $k \leq \frac{2}{3-\tau}$, there exist $\delta > 0$ such that
    \begin{equation*}
        \lim_{n \to \infty} \Prob{\frac{N^{(\star)}(H)}{n} \leq \delta} = 0.
    \end{equation*}
    \item When $k \geq \frac{2}{3-\tau}$, there exist $\delta > 0$ such that
    \begin{equation*}
        \lim_{n \to \infty} \Prob{\frac{N^{(\star)}(H)}{n^{k(3-\tau)/2}} \leq \delta} = 0.
    \end{equation*}
\end{itemize}
\end{corollary}

\begin{proof}
Observe that $\Nstar(H) \geq \Nstar(H,\Mab)$ for any $(\balpha,\bbeta) \in \mathcal{F}$. Moreover, $(\balpha,\bbeta) \equiv (1/2,0)$ and $(\balpha,\bbeta) \equiv (0,-1/d)$ are always in the feasible region for both optimization problems~\eqref{eq:optproblem} and~\eqref{eq:optproblem}--\eqref{eq:optproblem_extra_constraint}. Computing the exponent~\eqref{eq:fh} yields:
\begin{align*}
    f_H(\balpha,\bbeta) & = (3-\tau)k/2, \quad \text{if $(\balpha,\bbeta) \equiv (1/2,0)$}\\
    f_H(\balpha,\bbeta) & = 1, \quad \text{if $(\balpha,\bbeta) \equiv (0,-1/d)$}
\end{align*}
By \eqref{eq:subgraphGIRGprobbound}, there exists a function $g_1$ such that $\Prob{\Nstar(H;\Mab) \leq g_1(\varepsilon)n^{f_H(\balpha,\bbeta)}} \to 0,$ as $n\to \infty$. Hence, we conclude that there exists $\delta > 0$ such that 
$$\lim_{n \to \infty}\Prob{\Nstar(H) \leq \delta n^{\max(1,(3-\tau)k/2)}} = 0.$$
\end{proof}

Next, we derive an upper and a lower bound for the expected number of copies of $H$ in $\GIRGn$:

\begin{theorem}\label{thm:totalsubgraphcounts}
    Let $H$ be a subgraph of fixed size, let $\star \in \{\textup{ind},\textup{sub}\}$ and let $f_H^*$ be the optimal value of the optimization problem~\eqref{eq:optproblem} (if $\star = \textup{sub}$) or~\eqref{eq:optproblem}+\eqref{eq:optproblem_extra_constraint} (if $\star = \textup{ind}$).
    Then
    \begin{equation}\label{eq:ExpNHlowerbound}
    \liminf_{n \to \infty} \frac{\E{N^{(\star)}(H)}}{n^{f_H^*}} > 0.
    \end{equation}
    Moreover, conditionally on $\mathcal{E}_{\bar{\varepsilon}}$ defined in \eqref{eq:cutoff}, for any $\delta > 0$,
    \begin{equation}\label{eq:ExpNHupperbound}
    \limsup_{n \to \infty} \frac{\E{N^{(\star)}(H)\mid\mathcal{E}_{\bar{\varepsilon}}}}{n^{f_H^*+\delta}} < \infty.
    \end{equation}  
\end{theorem}
In other words, Theorem \ref{thm:totalsubgraphcounts} states that the expected subgraph count in $\GIRGn$ scales polynomially in the graph size, where the exponent is the solution of  \eqref{eq:optproblem}. Observe that from \eqref{eq:ExpNHupperbound} it also follows that, for any $\delta>0$, 
\begin{equation*}
    \lim_{n \to \infty} \Prob{N^{(\star)}(H) > n^{f_H^*+\delta}} = 0.
\end{equation*}
from direct application of Markov inequality.

\subsection{Results on particular subgraph classes}\label{sec:particularsubgraphs}

While the previous results hold in general for any  $H$ of fixed size, we can prove sharper convergence for specific classes of subgraphs.

First we focus on the case when $H$ is a tree, in which case the number of copies of $H$ as general subgraph in GIRGs shows the same asymptotic behaviour as in non-geometric scale free random graphs. In particular, it scales the same as in the inhomogeneous random graph (IRG) where the connection probability \eqref{eq:edgeprob} is replaced by $\tilde{p}_{uv} = \Theta\left( 1 \wedge \frac{w_{u} w_{v}}{\mu n}\right)$. We denote by $\widetilde{N}^{(\text{sub})}(H)$ the number of copies of $H$ in $\text{IRG}^{(n)}$ as general subgraph, for which the behavior of subgraph counts have been estimated in \cite{vdhofstad2021}.

\begin{theorem}\label{thm:subgraphtree}
If $H$ is a tree, then there exist $m,M>0$ such that
\begin{equation*}
    \lim_{n \to \infty} \Prob{m \leq \frac{\widetilde{N}^{(\text{sub})}(H)}{\Nsub(H)} \leq M} = 1.
\end{equation*}
\end{theorem}

The intuition behind this result is that the spatial embedding does not significantly influence the occurrence of pattern not containing cycles (i.e., trees). On the contrary, subgraphs that contain a Hamiltonian cycle undergo a geometric phase transition. Let $\Hk = ([k],\mathcal{E}_{\Hk})$ be the Hamiltonian cycle with edge set $\mathcal{E}_{\Hk}=\{(i,i+1)\}_{i \in [k]}$ using the convention that vertex $k+1$ coincides with vertex $1$. We say that $H$ is a Hamiltonian pattern, if there exist an automorphism of $H$ that contains $\Hk$.

\begin{theorem}\label{thm:hamiltonian}
    Let $H$ be a Hamiltonian pattern of size $k$. Let $\star \in \{\text{ind},\text{sub}\}$. Then there exist constants $I_{\text{G}}^{(\star)}(H),I_{\text{N}}^{(\star)}(H) \in (0,\infty)$ depending on the parameters of the GIRG such that:
    \begin{enumerate}[(i)]
        \item When $k < 1 + \frac{1}{3-\tau}$, the optimization problem for $H$ is solved by $\balpha \equiv 0, \bbeta \equiv -\frac{1}{d}$, and
        \begin{align}
            \frac{N^{(\star)}(H)}{n} \plim I_{\text{G}}^{(\star)}(H).
        \end{align}
        \item When $k > \frac{2}{3-\tau}$ and $k$ is odd, the optimization problem for $H$ is solved by $\balpha \equiv 1/2, \bbeta \equiv 0$, and
        \begin{align}\label{eq:oddcycle}
            \frac{N^{(\star)}(H)}{n^{k(3-\tau)/2}} \plim I_{\text{N}}^{(\star)}(H).
        \end{align}
    \end{enumerate} 
\end{theorem}

\subsection{MILP reformulation}\label{sec:resultMILP}

We now provide a strategy to easily solve the optimization problem~\eqref{eq:optproblem} or~\eqref{eq:optproblem}+\eqref{eq:optproblem_extra_constraint} by reformulating it as a mixed-integer linear optimization problem.
Furthermore, we will show that this formulation will allow to check whether a particular solution is unique, checking the assumption of Theorem~\ref{thm:subgraphGIRGprob}.

We first replace every minimum in~\eqref{eq:optproblem} by an auxiliary variable that is bounded from above by the two expressions in the minimum.
These variables are $\zeta_j$ for $\min_{i : i < j} \beta_{\{i,j\}}$ as well as $\delta_{\{i,j\}}$ for $\min(\alpha_i + \alpha_j - d \beta_{\{i,j\}} - 1, 0)$.
For the maximum in~\eqref{eq:optproblem_normal_constraint} we introduce a binary decision variable $z_{i,j,s}$ that indicates which of the constraints $\beta_{\{i,j\}} \leq \beta_{\{i,s\}}$ and $\beta_{\{i,j\}} \leq \beta_{\{j,s\}}$ must hold.
If that variable equals $1$, then~\eqref{eq:milpsubgraph:max1} and~\eqref{eq:milpsubgraph:max2} simplify to $\beta_{\{i,j\}} \leq \beta_{\{i,s\}}$ and $\beta_{\{i,j\}} \leq \beta_{\{j,s\}} + 1/d$, where the latter is redundant due to the bounds on the $\beta$-variables.
Similarly, if $z_{i,j,s} = 0$ then~\eqref{eq:milpsubgraph:max1} and~\eqref{eq:milpsubgraph:max2} simplify to $\beta_{\{i,j\}} \leq \beta_{\{i,s\}} + 1/d$ and $\beta_{\{i,j\}} \leq \beta_{\{j,s\}}$, where the former is redundant.
The resulting MILP reads as follows.
\begin{subequations}
  \label{eq:milpsubgraph}
    \begin{alignat}{7}
      & \text{max } ~\mathrlap{ (1-\tau) \sum_{i=1}^{k} \alpha_i  + d \sum_{j=1}^{k} \zeta_j + \gamma \sum_{ \{i,j\}\in \mathcal{E} } \delta_{\{i,j\}} } \label{eq:milp_obj}\\
    & \text{over } 
      & \alpha_i    &\in [0,1/(\tau-1)] 
        &\quad& \forall i \in [k] \label{eq:milpsubgraph:firstcons} \\
    & & \beta_{\{i,j\}}  &\in [-1/d,0] 
        &\quad& \forall i,j \in [k] : i < j \\
    & & \zeta_j     &\in [-1/d,0] 
        &\quad& \forall j \in \{2,3,\dotsc,k\} \\
    & & \delta_{\{i,j\}} &\in [-1,0] 
        &\quad& \forall \{i,j\} \in \mathcal{E} \\
    & & z_{i,j,s}     &\in \{0,1\} 
        &\quad& \forall i,j,s \in [k] : i < j,~ s \neq i,~ s \neq j \\
    & \text{s.t. }
      & \zeta_{j}   &\leq \beta_{\{i,j\}}
        &\quad& \forall i,j \in [k] : i < j \\
    & & \beta_{\{i,j\}}  &\leq \beta_{\{i,s\}} + (1 - z_{i,j,s}) / d 
        &\quad& \forall i,j,s \in [k] :~ i < j,~ s \neq i,~ s \neq j \label{eq:milpsubgraph:max1} \\
    & & \beta_{\{i,j\}}  &\leq \beta_{\{j,s\}} + z_{i,j,s} / d 
        &\quad& \forall i,j,s \in [k] :~ i < j,~ s \neq i,~ s \neq j \label{eq:milpsubgraph:max2} \\
    & & \delta_{\{i,j\}} &\leq \alpha_i + \alpha_j - d \beta_{\{i,j\}} - 1 
        &\quad& \forall \{i,j\} \in \mathcal{E} \label{eq:milpsubgraph:lastcons}
  \end{alignat}
\end{subequations}

Theorems~\ref{thm:subgraphGIRG}--\ref{thm:subgraphGIRGprob} show that for subgraphs for which the optimizer of~\eqref{eq:milpsubgraph} (resp.~\eqref{eq:milpsubgraph}+\eqref{eq:optproblem_extra_constraint}) is unique in terms of its $(\balpha,\bbeta)$ variables, there is a specific type of dominant subgraph structure.
That is, a copy of $H$ forms most likely on vertices of weights $\balpha$ and inter-distances $\bbeta$.

Hence, we are not just interested in an optimum solution of~\eqref{eq:milpsubgraph} or \eqref{eq:milpsubgraph}+\eqref{eq:optproblem_extra_constraint}, but also wish to know if the determined optimum solution is unique.
For linear optimization problems, sensitivity analysis can be used to assert uniqueness.
However, here we deal with a MILP and we only care about the uniqueness after projecting the optimal solution onto the $\balpha$- and $\bbeta$-variables. Therefore, we use the following approach.

First, we solve~\eqref{eq:milpsubgraph} or~\eqref{eq:milpsubgraph}+\eqref{eq:optproblem_extra_constraint} and obtain the optimizer $(\balpha^\star,\bbeta^\star)$ with objective value $z^\star$, restricting to only optimal solutions of the MILP.
Second, we add the constraint
\begin{equation}
  (1-\tau) \sum_{i=1}^{k} \alpha_i  + d \sum_{j=1}^{k} \zeta_j + \gamma \sum_{ \{i,j\}\in \mathcal{E} } \delta_{\{i,j\}} = z^\star
  \label{eq:milpobjectivecut}
\end{equation}
to the MILP.
Third, for every $i \in [k]$ and each $\sigma \in \{ \pm 1 \}$ we maximize $\sigma \alpha_i$ subject to~\eqref{eq:milpsubgraph:firstcons}--\eqref{eq:milpsubgraph:lastcons} and~\eqref{eq:milpobjectivecut} by replacing~\eqref{eq:milp_obj} by the new objective~$\sigma \in \{ \pm 1 \}$.
Fourth, for every $i,j \in [k] : i < j$ and each $\sigma \in \{ \pm 1 \}$ we maximize $\sigma \beta_{\{i,j\}}$ subject to~\eqref{eq:milpsubgraph:firstcons}--\eqref{eq:milpsubgraph:lastcons} and~\eqref{eq:milpobjectivecut}.
Constraint~\eqref{eq:milpobjectivecut} ensures that every feasible solution to any of these MILPs is optimal for~\eqref{eq:milpsubgraph} (resp.~\eqref{eq:milpsubgraph}+\eqref{eq:optproblem_extra_constraint}). Finally, we conclude that $(\balpha^\star,\bbeta^\star)$ is unique if and only if $(\balpha^\star,\bbeta^\star)$ is also the optimizer of all modified optimization problems.

If any such solution $(\balpha',\bbeta')$ is different from $(\balpha^\star,\bbeta^\star)$ then the optimal solution of~\eqref{eq:milpsubgraph} is not unique. Furthermore, if there exists an optimal solution of~\eqref{eq:milpsubgraph} different from $(\balpha^\star,\bbeta^\star)$ then it differs in at least one coordinate, where it is either smaller or larger.
Such a solution would constitute a feasible solution to each 
of these modified optimization problems. Furthermore, the difference in at least one coordinate implies that at least one $i$ or $i,j$ has a higher value of $\sigma \alpha_i$ or $\sigma \beta_{\{i,j\}}$ than $(\balpha^\star,\bbeta^\star)$.
Consequently, in at least one of these modified problems $(\balpha^\star,\bbeta^\star)$ will not be optimal with respect to the modified objective.
We conclude that the described approach correctly asserts uniqueness of $(\balpha^\star,\bbeta^\star)$ or outputs at least one alternative optimal solution.

\subsection{Numerical results}
We implemented~\eqref{eq:milpsubgraph} and~\eqref{eq:milpsubgraph}+\eqref{eq:optproblem_extra_constraint} as well as the approach for testing uniqueness using the open-source MILP-solver SCIP~\cite{SCIP9} through its Python interface~\cite{PYSCIPOPT}.

These results paint an interesting picture of the typical occurrences of all subgraphs of sizes 4 and 5. For example, for small $\tau$, dense subgraphs are often formed by $\sqrt{n}$-degree nodes that are far apart, while for larger values of $\tau$, dense subgraphs are often formed by low-degree vertices that are close in location. For less dense subgraphs, the picture is very subgraph-dependent, with some hubs and some lower-degree vertices, and also differences in the typical distances. Also in these sparser subgraphs, a phase transition in $\tau$ seems to appear, but the change in distances and weights of the involved vertices is more complex.

\begin{figure}[htpb]
  \stackunder[10pt]{\begin{tikzpicture}
        \begin{axis}[
        hide axis,
        scale only axis,
        height = 0pt,
        width = 0pt,
        point meta min = 0,
        point meta max = 1,
        legend cell align = {left},
        legend style = {draw = none},
        ]
        \addlegendimage{line width=4pt, alpha-1} \addlegendentry{$\alpha_v = 0$};
        \addlegendimage{line width=4pt, alpha-2} \addlegendentry{$\alpha_v = \frac{ \tau - 2 }{ \tau-1 }$};
        \addlegendimage{line width=4pt, alpha-3} \addlegendentry{$\alpha_v = 1/2$};
        \addlegendimage{line width=4pt, alpha-4} \addlegendentry{$\alpha_v = \frac{ 1 }{ \tau-1 }$};
        \addplot [draw=none] coordinates {(0,0)};
        \end{axis}
    \end{tikzpicture}
  }{Node colors ($\alpha$-values)}
    \hfill
  \stackunder[10pt]{\begin{tikzpicture}
        \begin{axis}[
        hide axis,
        scale only axis,
        height=0pt,
        width=0pt,
        colorbar horizontal,
        point meta min=-1.0,
        point meta max=0,
        colorbar style={width=0.68\textwidth}
        ]
        \addplot [draw=none] coordinates {(0,0)};
        \end{axis}
    \end{tikzpicture}
  }{Edge colors ($\beta$-values)}
  \caption{%
    Legend for other figures. Coloring of nodes and edges according to $\alpha$- and $\beta$-values, respectively.
    For subgraph plots, dashed edges indicate edges that are \emph{not} in the graph, while for induced subgraphs plots, these are drawn as dotted edges.
    In the figures, gray indicates that there also exist different optimal solutions.
  }
  \label{fig_legend}
\end{figure}
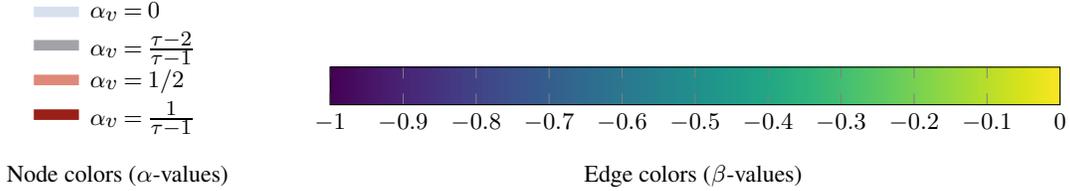

\vspace*{-7mm}
\setlength{\abovecaptionskip}{1pt plus 3pt minus 2pt}
\begin{figure}[htpb]
  \subfloat[
    (star)
    \label{fig_k4_taulow_subgraph_1}
  ]{
\begin{minipage}{0.135\textwidth} \begin{center} \begin{tikzpicture}
   \graphNode{v0}{45.0}{1} \graphNode{v1}{135.0}{1} \graphNode{v2}{225.0}{1} \graphNode{v3}{315.0}{4}
   \graphNonEdge{v0}{v1}{909.090909090909} \graphNonEdge{v0}{v2}{909.090909090909} \graphNonEdge{v1}{v2}{909.090909090909} \graphEdge{v0}{v3}{909.090909090909} \graphEdge{v1}{v3}{909.090909090909} \graphEdge{v2}{v3}{909.090909090909}
\end{tikzpicture} \end{center} \end{minipage}
  }%
  \hfill
  \subfloat[
    (path)
    \label{fig_k4_taulow_subgraph_2}
  ]{
\begin{minipage}{0.135\textwidth} \begin{center} \begin{tikzpicture}[nonUnique]
   \graphNode{v0}{45.0}{2} \graphNode{v1}{135.0}{4} \graphNode{v2}{225.0}{1} \graphNode{v3}{315.0}{1}
   \graphNonEdge{v0}{v2}{1000.0} \graphNonEdge{v1}{v3}{1000.0} \graphNonEdge{v2}{v3}{1000.0} \graphEdge{v0}{v1}{1000.0} \graphEdge{v0}{v3}{90.90909090909093} \graphEdge{v1}{v2}{909.090909090909}
\end{tikzpicture} \end{center} \end{minipage}
  }%
  \hfill
  \subfloat[
    (paw)
    \label{fig_k4_taulow_subgraph_3}
  ]{
\begin{minipage}{0.135\textwidth} \begin{center} \begin{tikzpicture}
   \graphNode{v0}{45.0}{1} \graphNode{v1}{135.0}{2} \graphNode{v2}{225.0}{2} \graphNode{v3}{315.0}{4}
   \graphNonEdge{v0}{v1}{1000.0} \graphNonEdge{v0}{v2}{1000.0} \graphEdge{v0}{v3}{909.090909090909} \graphEdge{v1}{v2}{181.81818181818187} \graphEdge{v1}{v3}{1000.0} \graphEdge{v2}{v3}{1000.0}
\end{tikzpicture} \end{center} \end{minipage}
  }%
  \hfill
  \subfloat[
    ($K_4$)
    \label{fig_k4_taulow_subgraph_4}
  ]{
\begin{minipage}{0.135\textwidth} \begin{center} \begin{tikzpicture}
   \graphNode{v0}{45.0}{3} \graphNode{v1}{135.0}{3} \graphNode{v2}{225.0}{3} \graphNode{v3}{315.0}{3}
   \graphEdge{v0}{v1}{1000.0} \graphEdge{v0}{v2}{1000.0} \graphEdge{v0}{v3}{1000.0} \graphEdge{v1}{v2}{1000.0} \graphEdge{v1}{v3}{1000.0} \graphEdge{v2}{v3}{1000.0}
\end{tikzpicture} \end{center} \end{minipage}
  }%
  \hfill
  \subfloat[
    (diamond)
    \label{fig_k4_taulow_subgraph_5}
  ]{
\begin{minipage}{0.135\textwidth} \begin{center} \begin{tikzpicture}[nonUnique]
   \graphNode{v0}{45.0}{4} \graphNode{v1}{135.0}{2} \graphNode{v2}{225.0}{4} \graphNode{v3}{315.0}{2}
   \graphNonEdge{v1}{v3}{1000.0} \graphEdge{v0}{v1}{1000.0} \graphEdge{v0}{v2}{1000.0} \graphEdge{v0}{v3}{1000.0} \graphEdge{v1}{v2}{1000.0} \graphEdge{v2}{v3}{1000.0}
\end{tikzpicture} \end{center} \end{minipage}
  }%
  \hfill
  \subfloat[
    (cycle)
    \label{fig_k4_taulow_subgraph_6}
  ]{
\begin{minipage}{0.135\textwidth} \begin{center} \begin{tikzpicture}[nonUnique]
   \graphNode{v0}{45.0}{4} \graphNode{v1}{135.0}{2} \graphNode{v2}{225.0}{4} \graphNode{v3}{315.0}{2}
   \graphNonEdge{v0}{v2}{1000.0} \graphNonEdge{v1}{v3}{1000.0} \graphEdge{v0}{v1}{1000.0} \graphEdge{v0}{v3}{1000.0} \graphEdge{v1}{v2}{1000.0} \graphEdge{v2}{v3}{1000.0}
\end{tikzpicture} \end{center} \end{minipage}
  }%
  \caption{%
    Optimal \textbf{subgraph} structures of all subgraphs with four nodes for \textbf{small} $\tau \in (2,2.5)$, $d=1$ and $\gamma = 2$.
    The legend is depicted in Figure~\ref{fig_legend}.
  }
  \label{fig_k4_taulow_subgraph}
\end{figure}
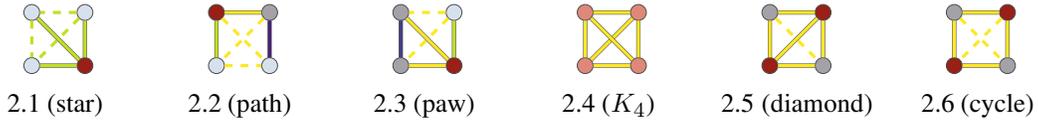
\vspace*{-7mm}

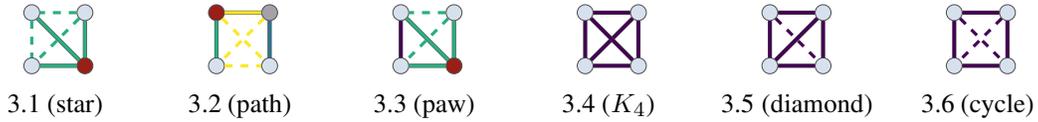
\begin{figure}[htpb]
  \subfloat[
    (star)
    \label{fig_k4_tauhigh_subgraph_1}
  ]{
\begin{minipage}{0.135\textwidth} \begin{center} \begin{tikzpicture}
   \graphNode{v0}{45.0}{1} \graphNode{v1}{135.0}{1} \graphNode{v2}{225.0}{1} \graphNode{v3}{315.0}{4}
   \graphNonEdge{v0}{v1}{625.0} \graphNonEdge{v0}{v2}{625.0} \graphNonEdge{v1}{v2}{625.0} \graphEdge{v0}{v3}{625.0} \graphEdge{v1}{v3}{625.0} \graphEdge{v2}{v3}{625.0}
\end{tikzpicture} \end{center} \end{minipage}
  }%
  \hfill
  \subfloat[
    (path)
    \label{fig_k4_tauhigh_subgraph_2}
  ]{
\begin{minipage}{0.135\textwidth} \begin{center} \begin{tikzpicture}[nonUnique]
   \graphNode{v0}{45.0}{2} \graphNode{v1}{135.0}{4} \graphNode{v2}{225.0}{1} \graphNode{v3}{315.0}{1}
   \graphNonEdge{v0}{v2}{1000.0} \graphNonEdge{v1}{v3}{1000.0} \graphNonEdge{v2}{v3}{1000.0} \graphEdge{v0}{v1}{1000.0} \graphEdge{v0}{v3}{375.0} \graphEdge{v1}{v2}{625.0}
\end{tikzpicture} \end{center} \end{minipage}
  }%
  \hfill
  \subfloat[
    (paw)
    \label{fig_k4_tauhigh_subgraph_3}
  ]{
\begin{minipage}{0.135\textwidth} \begin{center} \begin{tikzpicture}
   \graphNode{v0}{45.0}{1} \graphNode{v1}{135.0}{1} \graphNode{v2}{225.0}{1} \graphNode{v3}{315.0}{4}
   \graphNonEdge{v0}{v1}{624.9999999999999} \graphNonEdge{v0}{v2}{624.9999999999999} \graphEdge{v0}{v3}{625.0} \graphEdge{v1}{v2}{1.1102230246251565e-12} \graphEdge{v1}{v3}{624.9999999999999} \graphEdge{v2}{v3}{624.9999999999999}
\end{tikzpicture} \end{center} \end{minipage}
  }%
  \hfill
  \subfloat[
    ($K_4$)
    \label{fig_k4_tauhigh_subgraph_4}
  ]{
\begin{minipage}{0.135\textwidth} \begin{center} \begin{tikzpicture}
   \graphNode{v0}{45.0}{1} \graphNode{v1}{135.0}{1} \graphNode{v2}{225.0}{1} \graphNode{v3}{315.0}{1}
   \graphEdge{v0}{v1}{0.0} \graphEdge{v0}{v2}{0.0} \graphEdge{v0}{v3}{0.0} \graphEdge{v1}{v2}{0.0} \graphEdge{v1}{v3}{0.0} \graphEdge{v2}{v3}{0.0}
\end{tikzpicture} \end{center} \end{minipage}
  }%
  \hfill
  \subfloat[
    (diamond)
    \label{fig_k4_tauhigh_subgraph_5}
  ]{
\begin{minipage}{0.135\textwidth} \begin{center} \begin{tikzpicture}
   \graphNode{v0}{45.0}{1} \graphNode{v1}{135.0}{1} \graphNode{v2}{225.0}{1} \graphNode{v3}{315.0}{1}
   \graphNonEdge{v1}{v3}{0.0} \graphEdge{v0}{v1}{0.0} \graphEdge{v0}{v2}{0.0} \graphEdge{v0}{v3}{0.0} \graphEdge{v1}{v2}{0.0} \graphEdge{v2}{v3}{0.0}
\end{tikzpicture} \end{center} \end{minipage}
  }%
  \hfill
  \subfloat[
    (cycle)
    \label{fig_k4_tauhigh_subgraph_6}
  ]{
\begin{minipage}{0.135\textwidth} \begin{center} \begin{tikzpicture}
   \graphNode{v0}{45.0}{1} \graphNode{v1}{135.0}{1} \graphNode{v2}{225.0}{1} \graphNode{v3}{315.0}{1}
   \graphNonEdge{v0}{v2}{0.0} \graphNonEdge{v1}{v3}{0.0} \graphEdge{v0}{v1}{0.0} \graphEdge{v0}{v3}{0.0} \graphEdge{v1}{v2}{0.0} \graphEdge{v2}{v3}{0.0}
\end{tikzpicture} \end{center} \end{minipage}
  }%
  \caption{%
    Optimal \textbf{subgraph} structures of all subgraphs with four nodes for \textbf{large} $\tau \in (2.5,3)$, $d=1$ and $\gamma = 2$.
    The legend is depicted in Figure~\ref{fig_legend}.
  }%
  \label{fig_k4_tauhigh_subgraph}
\end{figure}
\vspace*{-7mm}

\begin{figure}[htpb]
  \subfloat[
    (star)
    \label{fig_k4_taulow_induced_1}
  ]{
\begin{minipage}{0.135\textwidth} \begin{center} \begin{tikzpicture}
   \graphNode{v0}{45.0}{1} \graphNode{v1}{135.0}{1} \graphNode{v2}{225.0}{1} \graphNode{v3}{315.0}{4}
   \graphNonEdgeInduced{v0}{v1}{909.090909090909} \graphNonEdgeInduced{v0}{v2}{909.090909090909} \graphNonEdgeInduced{v1}{v2}{909.090909090909} \graphEdge{v0}{v3}{909.090909090909} \graphEdge{v1}{v3}{909.090909090909} \graphEdge{v2}{v3}{909.090909090909}
\end{tikzpicture} \end{center} \end{minipage}
  }%
  \hfill
  \subfloat[
    (path)
    \label{fig_k4_taulow_induced_2}
  ]{
\begin{minipage}{0.135\textwidth} \begin{center} \begin{tikzpicture}[nonUnique]
   \graphNode{v0}{45.0}{4} \graphNode{v1}{135.0}{2} \graphNode{v2}{225.0}{1} \graphNode{v3}{315.0}{1}
   \graphNonEdgeInduced{v0}{v2}{1000.0} \graphNonEdgeInduced{v1}{v3}{1000.0} \graphNonEdgeInduced{v2}{v3}{1000.0} \graphEdge{v0}{v1}{1000.0} \graphEdge{v0}{v3}{909.090909090909} \graphEdge{v1}{v2}{90.90909090909093}
\end{tikzpicture} \end{center} \end{minipage}
  }%
  \hfill
  \subfloat[
    (paw)
    \label{fig_k4_taulow_induced_3}
  ]{
\begin{minipage}{0.135\textwidth} \begin{center} \begin{tikzpicture}
   \graphNode{v0}{45.0}{1} \graphNode{v1}{135.0}{2} \graphNode{v2}{225.0}{2} \graphNode{v3}{315.0}{4}
   \graphNonEdgeInduced{v0}{v1}{1000.0} \graphNonEdgeInduced{v0}{v2}{1000.0} \graphEdge{v0}{v3}{909.090909090909} \graphEdge{v1}{v2}{181.81818181818187} \graphEdge{v1}{v3}{1000.0} \graphEdge{v2}{v3}{1000.0}
\end{tikzpicture} \end{center} \end{minipage}
  }%
  \hfill
  \subfloat[
    ($K_4$)
    \label{fig_k4_taulow_induced_4}
  ]{
\begin{minipage}{0.135\textwidth} \begin{center} \begin{tikzpicture}
   \graphNode{v0}{45.0}{3} \graphNode{v1}{135.0}{3} \graphNode{v2}{225.0}{3} \graphNode{v3}{315.0}{3}
   \graphEdge{v0}{v1}{1000.0} \graphEdge{v0}{v2}{1000.0} \graphEdge{v0}{v3}{1000.0} \graphEdge{v1}{v2}{1000.0} \graphEdge{v1}{v3}{1000.0} \graphEdge{v2}{v3}{1000.0}
\end{tikzpicture} \end{center} \end{minipage}
  }%
  \hfill
  \subfloat[
    (diamond)
    \label{fig_k4_taulow_induced_5}
  ]{
\begin{minipage}{0.135\textwidth} \begin{center} \begin{tikzpicture}[nonUnique]
   \graphNode{v0}{45.0}{4} \graphNode{v1}{135.0}{2} \graphNode{v2}{225.0}{4} \graphNode{v3}{315.0}{2}
   \graphNonEdgeInduced{v1}{v3}{1000.0} \graphEdge{v0}{v1}{1000.0} \graphEdge{v0}{v2}{1000.0} \graphEdge{v0}{v3}{1000.0} \graphEdge{v1}{v2}{1000.0} \graphEdge{v2}{v3}{1000.0}
\end{tikzpicture} \end{center} \end{minipage}
  }%
  \hfill
  \subfloat[
    (cycle)
    \label{fig_k4_taulow_induced_6}
  ]{
\begin{minipage}{0.135\textwidth} \begin{center} \begin{tikzpicture}
   \graphNode{v0}{45.0}{3} \graphNode{v1}{135.0}{3} \graphNode{v2}{225.0}{3} \graphNode{v3}{315.0}{3}
   \graphNonEdgeInduced{v0}{v2}{1000.0} \graphNonEdgeInduced{v1}{v3}{1000.0} \graphEdge{v0}{v1}{1000.0} \graphEdge{v0}{v3}{1000.0} \graphEdge{v1}{v2}{1000.0} \graphEdge{v2}{v3}{1000.0}
\end{tikzpicture} \end{center} \end{minipage}
  }%
  \\
  \caption{%
    Optimal \textbf{induced subgraph} structures of all subgraphs with four nodes for \textbf{small} $\tau \in (2,2.5)$, $d=1$ and $\gamma = 2$.
    The legend is depicted in Figure~\ref{fig_legend}.
  }%
  \label{fig_k4_taulow_induced}
\end{figure}
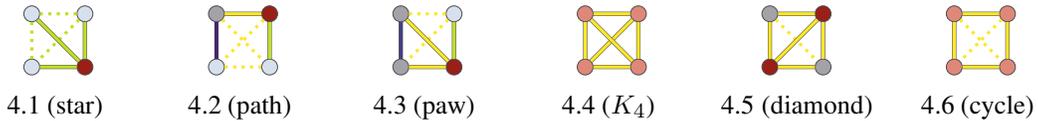
\vspace*{-3mm}

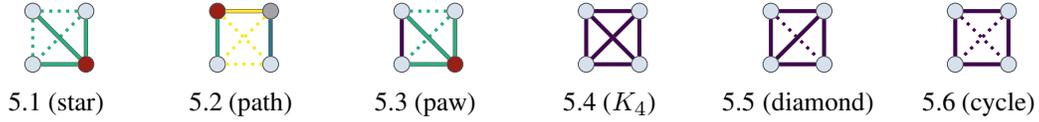
\begin{figure}[htpb]
  \subfloat[
    (star)
    \label{fig_k4_tauhigh_induced_1}
  ]{
\begin{minipage}{0.135\textwidth} \begin{center} \begin{tikzpicture}
   \graphNode{v0}{45.0}{1} \graphNode{v1}{135.0}{1} \graphNode{v2}{225.0}{1} \graphNode{v3}{315.0}{4}
   \graphNonEdgeInduced{v0}{v1}{625.0} \graphNonEdgeInduced{v0}{v2}{625.0} \graphNonEdgeInduced{v1}{v2}{625.0} \graphEdge{v0}{v3}{625.0} \graphEdge{v1}{v3}{625.0} \graphEdge{v2}{v3}{625.0}
\end{tikzpicture} \end{center} \end{minipage}
  }%
  \hfill
  \subfloat[
    (path)
    \label{fig_k4_tauhigh_induced_2}
  ]{
\begin{minipage}{0.135\textwidth} \begin{center} \begin{tikzpicture}[nonUnique]
   \graphNode{v0}{45.0}{2} \graphNode{v1}{135.0}{4} \graphNode{v2}{225.0}{1} \graphNode{v3}{315.0}{1}
   \graphNonEdgeInduced{v0}{v2}{1000.0} \graphNonEdgeInduced{v1}{v3}{1000.0} \graphNonEdgeInduced{v2}{v3}{1000.0} \graphEdge{v0}{v1}{1000.0} \graphEdge{v0}{v3}{375.0} \graphEdge{v1}{v2}{625.0}
\end{tikzpicture} \end{center} \end{minipage}
  }%
  \hfill
  \subfloat[
    (paw)
    \label{fig_k4_tauhigh_induced_3}
  ]{
\begin{minipage}{0.135\textwidth} \begin{center} \begin{tikzpicture}
   \graphNode{v0}{45.0}{1} \graphNode{v1}{135.0}{1} \graphNode{v2}{225.0}{1} \graphNode{v3}{315.0}{4}
   \graphNonEdgeInduced{v0}{v1}{625.0} \graphNonEdgeInduced{v0}{v2}{625.0} \graphEdge{v0}{v3}{625.0} \graphEdge{v1}{v2}{0.0} \graphEdge{v1}{v3}{625.0} \graphEdge{v2}{v3}{625.0}
\end{tikzpicture} \end{center} \end{minipage}
  }%
  \hfill
  \subfloat[
    ($K_4$)
    \label{fig_k4_tauhigh_induced_4}
  ]{
\begin{minipage}{0.135\textwidth} \begin{center} \begin{tikzpicture}
   \graphNode{v0}{45.0}{1} \graphNode{v1}{135.0}{1} \graphNode{v2}{225.0}{1} \graphNode{v3}{315.0}{1}
   \graphEdge{v0}{v1}{0.0} \graphEdge{v0}{v2}{0.0} \graphEdge{v0}{v3}{0.0} \graphEdge{v1}{v2}{0.0} \graphEdge{v1}{v3}{0.0} \graphEdge{v2}{v3}{0.0}
\end{tikzpicture} \end{center} \end{minipage}
  }%
  \hfill
  \subfloat[
    (diamond)
    \label{fig_k4_tauhigh_induced_5}
  ]{
\begin{minipage}{0.135\textwidth} \begin{center} \begin{tikzpicture}
   \graphNode{v0}{45.0}{1} \graphNode{v1}{135.0}{1} \graphNode{v2}{225.0}{1} \graphNode{v3}{315.0}{1}
   \graphNonEdgeInduced{v1}{v3}{0.0} \graphEdge{v0}{v1}{0.0} \graphEdge{v0}{v2}{0.0} \graphEdge{v0}{v3}{0.0} \graphEdge{v1}{v2}{0.0} \graphEdge{v2}{v3}{0.0}
\end{tikzpicture} \end{center} \end{minipage}
  }%
  \hfill
  \subfloat[
    (cycle)
    \label{fig_k4_tauhigh_induced_6}
  ]{
\begin{minipage}{0.135\textwidth} \begin{center} \begin{tikzpicture}
   \graphNode{v0}{45.0}{1} \graphNode{v1}{135.0}{1} \graphNode{v2}{225.0}{1} \graphNode{v3}{315.0}{1}
   \graphNonEdgeInduced{v0}{v2}{0.0} \graphNonEdgeInduced{v1}{v3}{0.0} \graphEdge{v0}{v1}{0.0} \graphEdge{v0}{v3}{0.0} \graphEdge{v1}{v2}{0.0} \graphEdge{v2}{v3}{0.0}
\end{tikzpicture} \end{center} \end{minipage}
  }%
  \caption{%
    Optimal \textbf{induced subgraph} structures of all subgraphs with four nodes for \textbf{large} $\tau \in (2.5,3)$, $d=1$ and $\gamma = 2$.
    The legend is depicted in Figure~\ref{fig_legend}.
  }%
  \label{fig_k4_tauhigh_induced}
\end{figure}
\vspace*{-7mm}


\setlength{\abovecaptionskip}{10pt plus 3pt minus 2pt}

\begin{figure}[htpb]
  \subfloat[
    (star)
    \label{fig_k5_taulow_subgraph_1}
  ]{
\begin{minipage}{0.135\textwidth} \begin{center} 
 \end{center} \end{minipage}
  }%
  \\

  \caption{%
    Optimal \textbf{subgraph} structures of all subgraphs with five nodes for \textbf{small} $\tau \in (2,2.6)$, $d=1$ and $\gamma = 2$.
    The legend is depicted in Figure~\ref{fig_legend}.
  }
  \label{fig_k5_taulow_subgraph}
\end{figure}

\begin{figure}[htpb]
  \subfloat[
    (star)
    \label{fig_k5_tauhigh_subgraph_1}
  ]{
\begin{minipage}{0.135\textwidth} \begin{center} %
 \end{center} \end{minipage}
  }%
  \\

  \caption{%
    Optimal \textbf{subgraph} structures of all subgraphs with five nodes for \textbf{large} $\tau \in (2.6,3)$, $d=1$ and $\gamma = 2$.
    The legend is depicted in Figure~\ref{fig_legend}.
  }
  \label{fig_k5_tauhigh_subgraph}
\end{figure}

\begin{figure}[htpb]
  \subfloat[
    (star)
    \label{fig_k5_taulow_induced_1}
  ]{
\begin{minipage}{0.135\textwidth} \begin{center} %
 \end{center} \end{minipage}
  }%
  \\

  \caption{%
    Optimal \textbf{induced subgraph} structures of all subgraphs with five nodes for \textbf{small} $\tau \in (2,2.6)$, $d=1$ and $\gamma = 2$.
    The legend is depicted in Figure~\ref{fig_legend}.
  }
  \label{fig_k5_taulow_induced}
\end{figure}

\begin{figure}[htpb]
  \subfloat[
    (star)
    \label{fig_k5_tauhigh_induced_1}
  ]{
\begin{minipage}{0.135\textwidth} \begin{center} %
 \end{center} \end{minipage}
  }%
  \\

  \caption{%
    Optimal \textbf{induced subgraph} structures of all subgraphs with five nodes for \textbf{large} $\tau \in (2.6,3)$, $d=1$ and $\gamma = 2$.
    The legend is depicted in Figure~\ref{fig_legend}.
  }
  \label{fig_k5_tauhigh_induced}
\end{figure}

\subsection{Discussion}\label{sec:discussion}

\paragraph*{Number of subgraphs outside the feasible region} 
In Section \ref{sec:optimalsubgraphs} we derived the number of subgraphs that appear on sets of the kind $\Mab$ in $\GIRGn$, where $(\balpha,\bbeta) \in \mathcal{F}$ is in the feasible region defined in~\eqref{eq:feasibleregion}. As already observed, whenever the variables $(\balpha,\bbeta)$ are chosen outside $\mathcal{F}$, the corresponding set $\Mab$ is empty with high probability. Thus, in this case the random variable $\Nstar(H;\Mab)$ converges to $0$ in probability as $n$ grows large. However, its expectation does not necessarily converge to 0 and, in fact, it could grow much larger than for any other feasible $(\balpha,\bbeta) \in \mathcal{F}$. One example is the wedge pattern $W = (\mathcal{V},\mathcal{E})$, with $\mathcal{V}=\{1,2,3\}, \mathcal{E} = \{(1,2),(2,3)\}$, for which it is possible to show that with high probability $\Nstar(H)$ is bounded by $o(n^{2/(\tau-1)+\delta})$ for arbitrary $\delta>0$, whereas $\E{\Nstar(H)}=\Omega(n^{4-\tau})$ due to the infinite-variance weights. This means that concentration results such as in Theorem~\ref{thm:subgraphGIRGprob} or Theorem~\ref{thm:hamiltonian} do not hold in general for all $H$. Therefore, investigating the limit of these non-concentrating subgraphs more closely could be interesting for further research.

\paragraph*{Phase transition on $\tau$}
Theorem~\ref{thm:hamiltonian} shows that the occurrence of any Hamiltonian pattern $H$ falls into one of two distinct geometric/non-geometric regimes, linking the pattern size and the power-law exponent. If $\tau$ is large, $H$ appears more frequently on low-degree vertices at small distances; when $\tau$ is small, then hubs of weight order $\sqrt{n}$ produce a large number of copies of $H$. Figures~\ref{fig_k4_taulow_subgraph},\ref{fig_k4_tauhigh_subgraph} and~\ref{fig_k5_taulow_subgraph}, \ref{fig_k5_tauhigh_subgraph} illustrate this phase transition for subgraphs of sizes $k=4$ and $5$. Often, the unique optimum is achieved at $(1/2,0)$ or $(0,-1/d)$ (depending on $\tau$) by the most dense subgraphs. When the maximum is achieved on $(1/2,0)$, it implies that, surprisingly, the way these subgraph form is not particularly affected by the geometry of the random graph. Theorem~\ref{thm:hamiltonian} contains a gap for the geometric/non-geometric transition between $(1/2,0)$ and $(0,-1/d)$ when $k\in[1+\frac{1}{3-\tau}, \frac{2}{3-\tau}]$. We conjecture that the phase transition is sharply separated by the curve $k = 2/(3-\tau)$, as this would match existing results for clique counts~\cite{michielan2021}.

\paragraph*{Values of $\alpha$ and $\beta$}
When $k \neq 2 / (3-\tau)$, that is outside the critical window, we conjecture that whenever there exist unique optimal $\alpha$-values of~\eqref{eq:optproblem} or~\eqref{eq:optproblem}+\eqref{eq:optproblem_extra_constraint} they take value in a finite set, $\alpha_i \in \{0, \frac{\tau-2}{\tau-1}, \frac{1}{2}, \frac{1}{\tau-1}\}$.
Experimentally, Figures~\ref{fig_k4_taulow_subgraph}--\ref{fig_k4_tauhigh_induced} and \ref{fig_k5_taulow_subgraph}--\ref{fig_k5_tauhigh_induced} also show that this is the case for subgraphs up to size 5 (results for subgraphs of size 6 are not included in this paper, but are also consistent.
The values $0$ and $1/(\tau-1)$ result from the domain constraint of alpha to be in $[0,1/(\tau-1)]$.
As~\eqref{eq:optproblem} optimizes a piecewise linear function, the optimizers over $\alpha$ are at the end points of the interval $[0, 1/(\tau-1)]$ or at the points where $\alpha_i+\alpha_j-d\beta_{ij}=1$.
This still leaves many possible options for $\alpha$, but we conjecture that the optimizers that are caused by these constraints only have $\alpha=1/2$ or $\alpha=(\tau-2)/(\tau-1)$.
Proving this is not immediate unfortunately, as the triangle inequality for the $\beta$ values makes it difficult to make small perturbations to optimal solutions.
Under this conjecture, the optimal $\beta$-values would also only be given by $0,-1/d$ (end points), and values where $\alpha_i+\alpha_j-d\beta_{ij}=1$.
For the conjectured values of $\alpha$, this would result in $d \cdot \beta_i \in \left\{ -1, \frac{2-\tau }{ \tau-1 }, - \frac{1}{2}, -\frac{1}{\tau-1}, \frac{\tau-3}{\tau-1}, \frac{\tau-3}{2 \tau-2 }, 0 \right\}$.

\begin{figure}[tpb]
  \subfloat[
    (path)
    \label{fig_k4_duplicates_2}
  ]{
\begin{minipage}{0.135\textwidth} \begin{center} \begin{tikzpicture}[nonUnique]
   \graphNode{v0}{45.0}{2} \graphNode{v1}{135.0}{4} \graphNode{v2}{225.0}{1} \graphNode{v3}{315.0}{1}
   \graphNonEdge{v0}{v2}{1000.0} \graphNonEdge{v1}{v3}{1000.0} \graphNonEdge{v2}{v3}{1000.0} \graphEdge{v0}{v1}{1000.0} \graphEdge{v0}{v3}{90.90909090909093} \graphEdge{v1}{v2}{909.090909090909}
\end{tikzpicture} \end{center} \end{minipage}
  }%
  \hfill
  \subfloat[
    (path)
    \label{fig_k4_duplicates_3}
  ]{
\begin{minipage}{0.135\textwidth} \begin{center} \begin{tikzpicture}[nonUnique]
   \graphNode{v0}{45.0}{4} \graphNode{v1}{135.0}{2} \graphNode{v2}{225.0}{1} \graphNode{v3}{315.0}{1}
   \graphNonEdge{v0}{v2}{1000.0} \graphNonEdge{v1}{v3}{1000.0} \graphNonEdge{v2}{v3}{1000.0} \graphEdge{v0}{v1}{1000.0} \graphEdge{v0}{v3}{909.090909090909} \graphEdge{v1}{v2}{90.90909090909093}
\end{tikzpicture} \end{center} \end{minipage}
  }%
  \hfill
  \subfloat[
    (diamond)
    \label{fig_k4_duplicates_6}
  ]{
\begin{minipage}{0.135\textwidth} \begin{center} \begin{tikzpicture}[nonUnique]
   \graphNode{v0}{45.0}{4} \graphNode{v1}{135.0}{2} \graphNode{v2}{225.0}{4} \graphNode{v3}{315.0}{2}
   \graphNonEdge{v1}{v3}{1000.0} \graphEdge{v0}{v1}{1000.0} \graphEdge{v0}{v2}{1000.0} \graphEdge{v0}{v3}{1000.0} \graphEdge{v1}{v2}{1000.0} \graphEdge{v2}{v3}{1000.0}
\end{tikzpicture} \end{center} \end{minipage}
  }%
  \hfill
  \subfloat[
    (diamond)
    \label{fig_k4_duplicates_7}
  ]{
\begin{minipage}{0.135\textwidth} \begin{center} \begin{tikzpicture}[nonUnique]
   \graphNode{v0}{45.0}{3} \graphNode{v1}{135.0}{3} \graphNode{v2}{225.0}{3} \graphNode{v3}{315.0}{3}
   \graphNonEdge{v1}{v3}{1000.0} \graphEdge{v0}{v1}{1000.0} \graphEdge{v0}{v2}{1000.0} \graphEdge{v0}{v3}{1000.0} \graphEdge{v1}{v2}{1000.0} \graphEdge{v2}{v3}{1000.0}
\end{tikzpicture} \end{center} \end{minipage}
  }%
  \hfill
  \subfloat[
    (cycle)
    \label{fig_k4_duplicates_8}
  ]{
\begin{minipage}{0.135\textwidth} \begin{center} \begin{tikzpicture}[nonUnique]
   \graphNode{v0}{45.0}{4} \graphNode{v1}{135.0}{2} \graphNode{v2}{225.0}{4} \graphNode{v3}{315.0}{2}
   \graphNonEdge{v0}{v2}{1000.0} \graphNonEdge{v1}{v3}{1000.0} \graphEdge{v0}{v1}{1000.0} \graphEdge{v0}{v3}{1000.0} \graphEdge{v1}{v2}{1000.0} \graphEdge{v2}{v3}{1000.0}
\end{tikzpicture} \end{center} \end{minipage}
  }%
  \hfill
  \subfloat[
    (cycle)
    \label{fig_k4_duplicates_9}
  ]{
\begin{minipage}{0.135\textwidth} \begin{center} \begin{tikzpicture}[nonUnique]
   \graphNode{v0}{45.0}{2} \graphNode{v1}{135.0}{4} \graphNode{v2}{225.0}{2} \graphNode{v3}{315.0}{4}
   \graphNonEdge{v0}{v2}{1000.0} \graphNonEdge{v1}{v3}{1000.0} \graphEdge{v0}{v1}{1000.0} \graphEdge{v0}{v3}{1000.0} \graphEdge{v1}{v2}{1000.0} \graphEdge{v2}{v3}{1000.0}
\end{tikzpicture} \end{center} \end{minipage}
  }%
  \\
  \caption{%
    \textbf{Several different optimal} \textbf{subgraph} structures of all subgraphs with four nodes for \textbf{small} $\tau \in (2,2.5)$, $d=1$ and $\gamma = 2$ for which the optimal solutions was not unique.
    The legend is depicted in Figure~\ref{fig_legend}.
  }%
  \label{fig_k4_different_optimal}
\end{figure}
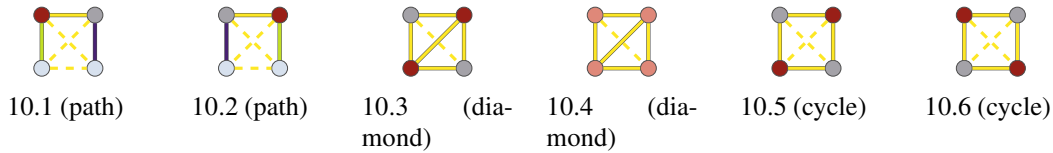
\paragraph*{Unique optimizers and predominant subgraphs}
The results of Corollary~\ref{thm:subgraphGIRG} and Theorem~\ref{thm:subgraphGIRGprob} hold for subgraphs with unique optimizers of~\eqref{eq:optproblem} or~\eqref{eq:optproblem}+\eqref{eq:optproblem_extra_constraint} is unique. This is the case for almost all subgraphs on 4 or 5 vertices. However, Figure~\ref{fig_k4_different_optimal} shows that non-unique optimizers are usually caused by symmetries in the subgraph. In that case, there is still a range of optimal subgraph structures, but rather than being defined on the vertices, they are defined on the edges. For example, on the cycle of size 4 for $\tau<5/2$, an optimal solution satisfies that the product of the vertex weights incident to any edge is of order $n$. Due to symmetry, this can arise in many different ways, and the extremal options are shown in Figure~\ref{fig_k4_different_optimal}, but in fact, the number of optimal solutions is infinite.

The reason that uniqueness of the optimizer is so important is the following. The optimal value of~\eqref{eq:optproblem} or~\eqref{eq:optproblem}+\eqref{eq:optproblem_extra_constraint} identifies the order of magnitude of subgraph counts within a specific region. Theorems~\ref{thm:subgraphGIRG}--\ref{thm:subgraphGIRGprob} then show that there exist optimal region(s), defined by the pair $(\balpha^*,\bbeta^*)$, maximizing the number of subgraphs in it. However, uniqueness of the optimizer does not immediately imply that the subgraphs in such maximal regions dominate all other regions combined. Practically, to demonstrate that subgraphs in the optimal region dominate all others, we need to show that for arbitrary $p_0 \in (0,1)$ there exists $\varepsilon > 0$ such that
\begin{equation}\label{eq:predominant}
    \lim_{n \to \infty}\Prob{\bv \in \M \mid H \subseteq \GIRGn|_{\bv}} \geq p_0.
\end{equation}
While this is not true in general, it is not difficult to show that~\eqref{eq:predominant} holds under the hypothesis of Theorem~\ref{thm:hamiltonian}, with a slight adaptation of its proof in Section \ref{sec:proofparticularsubgraphs}.

\paragraph*{Cycles of even size}
For Hamiltonian patterns of even size, Theorem~\ref{thm:hamiltonian}(ii) does not hold. The main reason is that for cycles of even size $k$ the solution to \eqref{eq:optproblem} is not unique. Indeed, in this case, two mirrored extreme optima appear, see Figure~\ref{fig_k4_different_optimal}. The convex combination of these two solution generates an infinite set of solutions. Introducing extra edges or considering the induced optimization problem may break the symmetry and result in a unique solution for the optimization problem. In this case, a more refined strategy is necessary to obtain similar results as in~\eqref{eq:oddcycle}. 

\paragraph*{Efficient subgraph counting algorithms}
A very rich line of research deals with the design of efficient algorithms to count or approximate subgraphs in networks \cite{kashtan2004,schreiber2005,ribeiro2022}. The knowledge of optimal set of vertices where subgraphs form more frequently may be used to improve significantly the efficiency of approximating algorithms. Indeed, one could use the following procedure: first, compute the (possibly unique) solution $(\balpha^*,\bbeta^*)$ of the general or induced optimization problem for a given pattern $H$; next, generate the induced graph over the set of vertices with degrees in $\cup_{\alpha \in \balpha^*} [\varepsilon n^{\alpha}, n^{\alpha}/\varepsilon]$, for some sensitivity parameter $\varepsilon \in (0,1)$; finally, run the algorithm for counting the occurrence of $H$ on the induced graph. As an example, consider the case when the optimum is uniquely achieved for $\alpha^* \equiv 1/2$. In this case, the size of the induced graph on vertices with order degree $\sim \sqrt{n}$ is $\Theta(n^{(3-\tau)/2})$; then, the number of edges in the induced graph is $O(n^{3-\tau})$, which is sub-linear with the size of the initial graph. Consequently, it is possible to define an approximated version of any node-based counting algorithm for $H$ improving the running time by $\Omega(n^{(\tau-1)/2})$, or of any edge-based counting algorithm improving the running time by $\Omega(n^{\tau-2})$. Moreover, the approximation error can be made arbitrarily small by tuning the parameter $\varepsilon$, at increased computational cost.

\section{Proofs for generic subgraphs}\label{sec:proofgenericsubgraphs}
Here we demonstrate results of Sections~\ref{sec:optimalsubgraphs}--\ref{sec:totalsubgraphs}. First we prove Proposition~\ref{thm:optproblem_gamma_independent}, showing that the parameter $\gamma$ does not play any role in the optimization problem~\ref{eq:optproblem}. Next, we prove the central Lemma~\ref{lemma:expNHMab}, showing that the subgraph counts over sets of the form $\Mab$ are particularly easy to compute.To demonstrate Theorem~\ref{thm:subgraphGIRGprob} we need to show that the subgraph counts are self-averaging random variables where we use a second moment method. Finally, we will prove upper and lower bounds for the total number of subgraphs in Theorem~\ref{thm:totalsubgraphcounts}.

\begin{proof}[Proof of Proposition~\ref{thm:optproblem_gamma_independent}]
Suppose that $\alpha_i + \alpha_j < 1 + d \beta_{\{i,j\}}$ for some $i,j$.
We will show that the objective function $f_H$ defined by~\eqref{eq:optproblem_objective} achieves a higher value when we decrease $\beta_{\{i,j\}}$ by a sufficiently small quantity $x>0$.
When decreasing $\beta_{\{i,j\}}$, we need to ensure that all constraints given by the triangular inequality are still satisfied.
To do so, we apply the following policy that keeps track of what $\beta_{ik}$ or $\beta_{jk}$ values are the same as $\beta_{ij}$, while one other distance in the triple $\{i,j,k\}$ is larger. This would imply that decreasing $\beta_{ij}$ can only be feasible when also decreasing $\beta_{ik}$ or $\beta_{jk}$:
\begin{enumerate}
    \item Start by declaring all edges $(a,b)$ with $a,b \in [k]$ inactive.
    \item Activate edge $(i,j)$ and consider all (unordered) triplets of distinct elements $(i,j,h)$ with $h \in [k]$:
    \begin{itemize}
        \item If $\beta_{\{i,j\}} \leq \beta_{\{i,h\}} = \beta_{\{j,h\}}$, keep $(i,h)$ and $(j,h)$ inactive.
        \item If $\beta_{\{i,h\}} < \beta_{\{i,j\}} = \beta_{\{j,h\}}$, switch $(j,h)$ to active, and keep $(i,h)$ inactive.
        \item If $\beta_{\{j,h\}} < \beta_{\{i,j\}} = \beta_{\{i,h\}}$, switch $(i,h)$ to active, and keep $(j,h)$ inactive.
    \end{itemize}
    \item For each newly activated edge $(i',j')$, repeat the operation in step 2, considering all triplets $(i',j',h)$.
\end{enumerate}
By construction, decreasing $\bbeta$ with a sufficiently small amount $x>0$ over all active edges is feasible with respect to the triangular inequality. 

W.l.o.g. assume $\alpha_1 + \alpha_2 < 1 + d \beta_{\{1,2\}}$. Decrease the values in $\bbeta$ over the activated edges by $x>0$ sufficiently small, accordingly to the activation policy described above, starting from edge $(1,2)$. Observe that, for each $t \geq 3$, there are exactly two possibilities:
\begin{itemize}
    \item  $\min_{s < t}\beta_{\{s,t\}} < \beta_{\{1,2\}}$. Let $s'$ be a vertex that achieves the minimum. Then, the policy keeps the edge $(s',t)$ inactive, as it only activates edges with $\beta$ value equal to $\beta_{\{1,2\}}$;
    \item $\min_{s < t}\beta_{\{s, t\}} \geq \beta_{\{1,2\}}$. Denote by $(s',t)$ the edge achieving the minimum. In this case, observe that vertex $t$ is at least at distance $\Omega(n^{\beta_{\{1,2\}}})$ from edge $(1,2)$. In turn, it follows from our activation policy that $t$ is at distance $\Omega(n^{\beta_{\{1,2\}}})$ from any other activated edge. In particular, this implies that the edge $(s',t)$ will never be activated.
\end{itemize}
In other words, whenever we decrease $\beta_{\{1,2\}}$, the triangular inequality constraints are satisfied and the quantities $\min_{s<t} \beta_{\{s,t\}}$ are kept fixed, when applying our decreasing policy.\\

Finally, observe that when decreasing $\beta_{\{1,2\}}$ (and eventually other $\beta_{\{s,t\}}$ for some $s,t \in [k]$) by $x>0$, the third sum in the objective function $f_H$ increases at least by $\gamma d x$. On the other hand, among the $k-1$ terms in the second sum of the objective function $f_H$, only the first (corresponding to index $j=2$) decreases by a quantity $d x$, whereas all the others remain unchanged. Therefore, decreasing $\beta_{\{1,2\}}$ increases the value of $f_H$ by $(\gamma-1)dx > 0$.
\end{proof}

\begin{proof}[Proof of Lemma \ref{lemma:expNHMab}]
    The definitions in \eqref{eq:Nha} yield
    \begin{equation*}
        \mathbb{E}[\Nind(H;\Mab)] = \sum_{\bv \in V_k} \mathbb{P}(H = \GIRG|_{\bv} \mid \bv \in \Mab) \mathbb{P}(\bv \in \Mab)
    \end{equation*}
    \begin{equation*}
        \mathbb{E}[\Nsub(H;\Mab)] = \sum_{\bv \in V_k} \mathbb{P}(H \subseteq \GIRG|_{\bv} \mid \bv \in \Mab) \mathbb{P}(\bv \in \Mab).
    \end{equation*}
    Given $\bv \in \Mab$, each vertex $v_i \in \bv$ has weight in $I_\varepsilon(n^{\alpha_i})$ with a probability which is proportional to $n^{(1-\tau)\alpha_i}$ due to~\eqref{eq:weightdistribution}, independently of all the other vertices in $\bv$. 
The probability that the inter-distance between any two selected vertices $v_i,v_j \in \bv$ falls in $I_\varepsilon(n^{\beta_{\{i,j\}}})$ is slightly more complicated. We calculate this probability starting from vertex $v_1$, which is w.l.o.g. given a fixed location on the torus. Next, consider $v_2$: the region in which the distance between $v_2$ and $v_1$ is in $I_\varepsilon(n^{\beta_{\{1,2\}}})$ is proportional to $n^{d\beta_{\{1,2\}}}$, as the positions are sampled uniformly. We then proceed to vertex $v_{3}$. The distance from $v_3$ to $v_1$ and $v_2$ should be inside $I_\varepsilon(n^{\beta_{\{1,3\}}})$ and $I_\varepsilon(n^{\beta_{\{2,3\}}})$, hence the area in which $v_3$ satisfies this is proportional to $\min(n^{d\beta_{\{1,3\}}},n^{d\beta_{\{2,3\}}})$. Iterating this procedure, combined with the fact that positions are sampled uniformly on the torus, yields that the probability that the vertex set $\bv\in \Mab$ is proportional to
\begin{equation}\label{eq:ProbvinMab}
    \mathbb{P}(\bv \in \Mab) \propto \prod_{i\in[k]}n^{(1-\tau)\alpha_i}\prod_{j\in[k]}\min_{i:i<j}n^{d\beta_{\{i,j\}}}.
\end{equation}

Furthermore, since $v_1,...,v_k\in \Mab$, the connection probability between $v_i$ and $v_j$ is proportional to $\min(n^{\alpha_i+\alpha_j - d \beta_{\{i,j\}} - 1}, 1)^{\gamma}$. This yields
\begin{equation}\label{eq:ProbHsubMab}
   \mathbb{P}(H \subseteq \GIRG|_{\bv} \mid \bv \in \Mab) \propto \prod_{i,j:(i,j) \in \mathcal{E}} \min(n^{\alpha_i+\alpha_j - d \beta_{\{i,j\}} - 1}, 1)^{\gamma}.
\end{equation}

Then, \eqref{eq:ProbvinMab} and \eqref{eq:ProbHsubMab}, and the fact that $\bv$ can be chosen in $n!/(n-k)!$ 
ways, gives
\begin{equation*}
    \mathbb{E}[\Nsub(H;\Mab)] \propto n^{f_H(\balpha,\bbeta)},
\end{equation*}
where
\begin{equation*}\label{eq:Mabexponent}
    f_H(\balpha,\bbeta) = k + (1-\tau)\sum_{i=1}^k \alpha_i + d\sum_{j=2}^k\min_{i:i<j}\beta_{\{i,j\}} + \gamma \sum_{i,j:(i,j) \in \mathcal{E}} \min(\alpha_i+\alpha_j - d \beta_{\{i,j\}} - 1,0).
\end{equation*}

For $\mathbb{E}[\Nind(H;\Mab)]$, we compute $\mathbb{P}(H = \GIRG|_{\bv} \mid \bv \in \Mab)$, similarly to \eqref{eq:ProbHsubMab} as
\begin{align}\label{eq:ProbHindMab}
   \mathbb{P}(H = \GIRG|_{\bv} \mid \bv \in \Mab) &\propto \prod_{(i,j) \in \mathcal{E}} \min(n^{\alpha_i+\alpha_j - d \beta_{\{i,j\}} - 1}, 1)^{\gamma} \prod_{(i,j) \not\in \mathcal{E}} \big( 1 - \min(n^{\alpha_i+\alpha_j - d \beta_{\{i,j\}} - 1}, 1)^{\gamma} \big) \nonumber\\
   &\propto \prod_{(i,j) \in \mathcal{E}} \min(n^{\alpha_i+\alpha_j - d \beta_{\{i,j\}} - 1}, 1)^{\gamma} \prod_{(i,j) \not\in \mathcal{E}} \max(1 - n^{\alpha_i+\alpha_j - d \beta_{\{i,j\}} - 1}, 0)^{\gamma} .\nonumber
\end{align}
In particular, whenever $\alpha_i+\alpha_j - d \beta_{\{i,j\}} - 1 > 0$ for any $(i,j) \not \in \mathcal{E}$, $\prob (H = \GIRG|_{\bv})=0$. Indeed, if a pair $(i,j) \not \in \mathcal{E}$ existed such that $\alpha_i+\alpha_j - d \beta_{\{i,j\}} - 1 > 0$, then the edge $(v_i,v_j)$ would be present with probability one, by~\eqref{eq:edgeprob}.
On the other hand, if $\alpha_i+\alpha_j - d \beta_{\{i,j\}} - 1 \leq 0$ for all $(i,j) \not \in \mathcal{E}$, then the expected number of copies of $H$ as induced subgraph is asymptotically equal to the expected number of copies of $H$ as general subgraph and
\begin{equation*}
    \mathbb{E}[\Nind(H;\Mab)] \propto
    \begin{cases}
        n^{f_H(\balpha,\bbeta)} & \text{ if $\alpha_i+\alpha_j - d \beta_{\{i,j\}} - 1 \leq 0, \forall (i,j) \not \in \mathcal{E}$}, \\
        0 & \text{ otherwise}
    \end{cases}
\end{equation*}
follows, which concludes the proof.
\end{proof}

Next, we prove that when the optimal solution is either $(\balpha^*,\bbeta^*) \equiv (1/2, 0)$ or $(\balpha^*,\bbeta^*) \equiv (0, -1/d)$, the convergence of expectations for the number of subgraphs in the optimal region can be sharpened to a convergence in probability. We first introduce an auxiliary lemma.

\begin{lemma}\label{lemma:selfaveraging}
Let $H$ be a graph, $\star \in \{\text{ind},\text{sub}\}$ and $(\balpha,\bbeta) \equiv (1/2, 0)$ or $(\balpha,\bbeta) \equiv (0, -1/d)$. Then, 
\begin{equation*}
    \lim_{n \to \infty} \frac{\Var{N^{(\star)}(H;\Mab)}}{\E{N^{(\star)}(H;\Mab)}^2} = 0.
\end{equation*}
\end{lemma}

\begin{proof}
For any fixed combination of $k$ vertices $\bv$ we define the events
\begin{gather*}
    A_{\bv} = \{H \subseteq \GIRGn|_{\bv}\},\quad
    B_{\bv} = \{\bv \in \Mab\},
\end{gather*}
so that $N^{\text{sub}}(H,\Mab) = \sum_{\bv \in V_k} \mathbbm{1}_{\{A_{\bv},B_{\bv}\}}$. Then, we write its variance as
\begin{equation*}
\begin{split}
    \Var{N^{\text{sub}}(H,\Mab)} &= \Var{\sum_{\bv \in V_k} \mathbbm{1}_{\{A_{\bv},B_{\bv}\}} }\\
    &= \sum_{\bv,\bu \in V_k} \Cov{ \mathbbm{1}_{\{A_{\bv},B_{\bv}\}},\mathbbm{1}_{\{A_{\bu},B_{\bu}\}}} \\
    &= \sum_{\bv,\bu \in V_k} \Prob{A_{\bv},B_{\bv},A_{\bu},B_{\bu}} - \Prob{A_{\bv},B_{\bv}}\Prob{(A_{\bu},B_{\bu}}.
\end{split}
\end{equation*}
Observe that if $\bv \cap \bu = \varnothing$, then the covariance between $\mathbbm{1}_{\{A_{\bv},B_{\bv}\}}$ and $\mathbbm{1}_{\{A_{\bu},B_{\bu}\}}$ is $0$. Therefore, we restrict to the case $|\bv \cap \bu| = s$, with $s \geq 1$. Furthermore, for all pairs $\bv, \bu$ we use the bound 
$$\Prob{A_{\bv},B_{\bv},A_{\bu},B_{\bu}} - \Prob{A_{\bv},B_{\bv})\Prob{A_{\bu},B_{\bu}}} \leq \Prob{A_{\bv},B_{\bv},A_{\bu},B_{\bu}} \leq \Prob{B_{\bv},B_{\bu}}.$$
Following the same computations as in the proof of Lemma \ref{lemma:expNHMab} leading to Equation \eqref{eq:ProbvinMab}, and assuming that $\bv$ and $\bu$ intersect in $s$ elements, we obtain
\begin{align*}
    \Prob{B_{\bv},B_{\bu}} &= \Theta(n^{k(1-\tau)\alpha + (k-1)d\beta}) \cdot \Theta(n^{(k-s)(1-\tau)\alpha + (k-s)d\beta}) \nonumber \\
    &= \Theta(n^{(2k-s)(1-\tau)\alpha + (2k-s-1)d\beta}),
\end{align*}
where $(\alpha,\beta)=(1/2,0)$ or $(\alpha,\beta)=(0, -1/d)$. 
Since there are in total $\Theta(n^{2k-s})$ ways to choose the vertices in $\bv$ and $\bu$, we obtain
\begin{align*}
    \sum_{\bv,\bu : |\bv \cap \bu| = s} \Cov{\mathbbm{1}_{\{A_{\bv},B_{\bv}\}},\mathbbm{1}_{\{A_{\bu},B_{\bu}\}}} &\leq \Theta(n^{2k-s}) \cdot \Theta(n^{(2k-s)(1-\tau)\alpha^* + (2k-s-1)d\beta^*}) \nonumber \\
    &= \Theta(n^{(2k-s) + (2k-s)(1-\tau)\alpha^* + (2k-s-1)d\beta^*}),
\end{align*}
and taking the sum over all possible values for $s$,
\begin{align*}
    \Var{N^{\text{sub}}(H;\Mab)} = O(\max_{s\in{1,\dots,k}}n^{(2k-s) + (2k-s)(1-\tau)\alpha + (2k-s-1)d\beta}).
\end{align*}
By Lemma \ref{lemma:expNHMab}, $\E{\Nsub(H;\Mab)} = \Omega(n^{k + k(1-\tau)\alpha + (k-1)d\beta})$. Therefore
\begin{align*}
    \frac{\Var{\Nsub(H;\Mab)}}{\E{\Nsub(H;\Mab)}^2} &= \frac{O(\max_{s\in{1,\dots,k}}n^{(2k-s) + (2k-s)(1-\tau)\alpha + (2k-s-1)d\beta})}{\Omega(n^{2k + 2k(1-\tau)\alpha + 2(k-1)d\beta})} \nonumber \\
    &= O(\max_{s\in{1,\dots,k}}n^{-s -s(1-\tau)\alpha + (1-s)d\beta}).
\end{align*}
Now, $g(s) = -s -s(1-\tau)\alpha + (1-s)d\beta$ is negative for all $s \geq 1$, as 
\begin{itemize}
    \item when $(\alpha,\beta) = (0,-1/d)$ then $g(s) = -s -(1-s) = -1$,
    \item when $(\alpha,\beta) = (1/2,0)$ then $g(s) = -s -s(1-\tau)/2 = -s(3-\tau)/2$.
\end{itemize}
This proves the claim. In the case of induced subgraph, the proof works by replacing $\Nsub$ with $\Nind$ and defining instead $A_{\bv}:= \{H = \GIRGn|_{\bv}\}$ at the beginning of the proof.
\end{proof}

\begin{proof}[Proof of Theorem \ref{thm:subgraphGIRGprob}]
The Chebyshev inequality states that for any random variable $X$ with finite variance and $\delta>0$, 
$$\Prob{\left|\frac{X}{\E{X}} - 1 \right| \geq \delta} \leq \frac{\Var{X}}{\delta^2\E{X}^2}.$$
Combining this with Lemma \ref{lemma:selfaveraging} yields that for any $\delta > 0$
$$\Prob{\left|\frac{\Nstar(H;\M)}{\E{\Nstar(H;\M)}} - 1 \right| \geq \delta} \longrightarrow 0, \quad \text{as $n\to \infty$}.$$
Combining this with Lemma \ref{lemma:expNHMab} proves \eqref{eq:subgraphGIRGprobbound}.

Next, let $(\balpha,\bbeta) \in \mathcal{F}\setminus \{(\balpha^*, \bbeta^*)\}$. Lemma \ref{lemma:expNHMab} shows that there exists a function $h_2$ such that
$$\frac{\E{\Nstar(H;\Mab)}}{n^{f_H(\balpha,\bbeta)}} \leq h_2(\varepsilon).$$
Then the Markov inequality yields
\begin{equation}
\begin{split}
    \lim_{n \to \infty}\Prob{\Nstar(H;\Mab) \geq \delta n^{f_H(\balpha^*,\bbeta^*)}} &\leq \lim_{n\to\infty}\frac{\E{\Nstar(H;\Mab)}}{\delta n^{f_H(\balpha^*,\bbeta^*)}} \\ &\leq \lim_{n\to\infty}\frac{h_2(\varepsilon)n^{f_H(\balpha,\bbeta)}}{\delta n^{f_H(\balpha^*,\bbeta^*)}} = 0,\label{eq:firstprobbound}
\end{split}
\end{equation}
for any $\delta > 0$. Moreover, by \eqref{eq:subgraphGIRGprobbound}, there exists a function $g_1$ such that
\begin{equation}\label{eq:secondprobbound}
    \Prob{\Nstar(H;\M) \geq g_1(\varepsilon) n^{f_H(\balpha^*,\bbeta^*)}} \longrightarrow 1, \quad \text{as $n \to \infty.$}
\end{equation}
Therefore, combining \eqref{eq:firstprobbound} and \eqref{eq:secondprobbound}, we obtain
\begin{equation*}
    \Prob{\frac{\Nstar(H;\Mab)}{\Nstar(H;\M)} \geq \delta} \longrightarrow 0, \quad \text{as $n \to \infty$},
\end{equation*}
which concludes the proof of the theorem.
\end{proof}

Finally, we demonstrate the result for the total subgraph counts. The lower bound follows immediately from Lemma~\ref{lemma:expNHMab}, whereas the upper bound involves building a proper covering of the feasible region $\mathcal{F}$.

\begin{proof}[Proof of Theorem \ref{thm:totalsubgraphcounts}]
\textbf{Lower bound.} The lower bound in \eqref{eq:ExpNHlowerbound} follows from the fact that, if $(\balpha^*,\bbeta^*)$ is a solution of the general or induced optimization problem, then
\begin{equation*}
    \E{\Nstar(H)} \geq \En{\Nstar(H;\M)} = \Omega(n^{f_H^*})
\end{equation*}
by Lemma \ref{lemma:expNHMab}.\\

\textbf{Upper bound.} To prove the upper bound in \eqref{eq:ExpNHupperbound}, we construct a covering for the feasible region. On $\mathcal{E}$, we have
\begin{equation*}
    \Nstar(H) = \Nstar\left(H; \bigcup_{(\balpha,\bbeta) \in \mathcal{F}} M_{\bar{\varepsilon}}^{(\balpha,\bbeta)} \right) \qquad w.h.p.
\end{equation*}
for some $\bar{\varepsilon}>0$, where $\mathcal{F}$ is the feasible region defined in \eqref{eq:feasibleregion}.

Covering $[n^a,n^b]$, for $a,b\in \mathbb{R}$ with intervals of the type $I_{\bar{\varepsilon}}(n^{\zeta})$, with $\zeta \in [a,b]$ translates to covering $[a\log(n),b\log(n)]$ in logarithmic scale with intervals of the type $[\zeta \log(n) + \log(\bar{\varepsilon}), \zeta \log(n) - \log(\bar{\varepsilon})]$ of length $-2 \log(\bar{\varepsilon})$. Therefore, such a covering has cardinality $\left\lceil \frac{(b-a) \log(n)}{-2 \log(\bar{\varepsilon})} \right\rceil$. In particular, vertex weights take values in the interval $[1,n^{1/(\tau-1)}]$, which can be covered with $\left\lceil \frac{\log(n)}{-2 (\tau-1) \log(\bar{\varepsilon})} \right\rceil$ sub-intervals of the form $I_{\bar{\varepsilon}}(n^{\zeta})$. Similarly, all inter-distances $||x_{v_i}-x_{v_j}||$ take value in the interval $[n^{-1/d},1]$ which is covered by $\left\lceil \frac{\log(n)}{-2 d \log(\bar{\varepsilon})} \right\rceil$ sub-intervals.
In other words, $\bigcup_{(\balpha,\bbeta) \in \mathcal{F}} M_{\bar{\varepsilon}}^{(\balpha,\bbeta)}$ admits a finite covering over a set $\mathcal{C} \subset \mathcal{F}$, with cardinality $$|\mathcal{C}| = \left\lceil \frac{\log(n)}{-2 (\tau-1) \log(\bar{\varepsilon})} \right\rceil^k \left\lceil \frac{\log(n)}{-2 d\log(\bar{\varepsilon})} \right\rceil^{\binom{k}{2}} = \Theta(\log(n)^{k(k+1)/2}).$$ Then,
\begin{equation*}
    \Nstar(H) \leq \sum_{(\balpha,\bbeta) \in \mathcal{C}} \Nstar(H;M_{\bar{\varepsilon}}^{(\balpha,\bbeta)} ), \quad w.h.p.
\end{equation*}
and from Lemma \ref{lemma:expNHMab} we conclude that
\begin{equation*}
    \En{\Nstar(H)} = \sum_{(\balpha,\bbeta) \in \mathcal{C}} \En{\Nstar (H;M_{\bar{\varepsilon}}^{(\balpha,\bbeta)})} \leq |\mathcal{C}| \cdot \Theta(n^{f_H^*}) = o(n^{f_H^* + \delta})
\end{equation*}
for any $\delta > 0$.
\end{proof}

\section{Proofs for specific subgraphs}\label{sec:proofparticularsubgraphs}
In this section we prove results of Section~\ref{sec:particularsubgraphs}. 

\subsection{Proof for trees}
We need to show Theorem~\ref{thm:subgraphtree}, that is, tree counts coincide asymptotically in $\GIRG$ and IRG models. In particular we aim to compare $\GIRGn$, where the connection probability between any pair of nodes $i,j$ is $p_{ij}=\min\left\{1, \frac{w_i w_j}{n \mu ||x_i-x_j||^d}\right\}^{\gamma}$, with IRG$_n$, where the connection probability is $\tilde{p}_{ij} = \Theta\left( 1 \wedge \frac{w_{v_i} w_{v_j}}{\mu n}\right)$. We start by proving that in trees the spatial edge dependency disappears.

\begin{lemma}\label{lemma:treeedgesindependence}
Let $T = (V_T,E_T)$ be a tree. Let $(X_i)_{i \in V_T}$ be a sequence of i.i.d. uniformly distributed random variables on the $d$-dimensional torus $\mathbb{T}^d$. For all $i,j \in V_T$ let $f_{ij}:[0,1/2]\to \mathbb{R}$ be a real function. Then,
\begin{equation}\label{eq:treeedgesindependence}
    \mathbb{E}\left[ \prod_{(i,j) \in E_T} f_{ij}(||X_i-X_j||) \right] = \prod_{(i,j) \in E_T} \mathbb{E}\left[f_{ij}(||X_i-X_j||) \right].
\end{equation}
\end{lemma}

\begin{proof}
The proof relies on the fact that, if $Y_1$ and $Y_2$ are uniformly distributed on $\mathbb{T}^d$ and $g$ is any real function, then 
\begin{equation*}
    \mathbb{E}[g(||Y_1 - Y_2||)] = \mathbb{E}[g(||Y_1||)].
\end{equation*}
Indeed, conditioning on $Y_2$ yields
\begin{align*}
    \mathbb{E}\left[\mathbb{E}[g(||Y_1 - Y_2||)|Y_2]\right] &= \int_{\mathbb{T}^d} \left( \int_{\mathbb{T}^d} g(||y_1 - y_2||) \dd y_1 \right) \dd y_2 \nonumber\\
    &= \int_{\mathbb{T}^d} \left( \int_{\mathbb{T}^d} g(||y_1||) \dd y_1 \right) \dd y_2 = \mathbb{E}[g(||Y_1||)],
\end{align*}
where the middle equality follows by definition of the infinity norm and symmetry of the uniform distribution.

We prove \eqref{eq:treeedgesindependence} by induction on the number of edges $|E_T|$. If $|E_T|=1$, the equation is trivial. Then, suppose that $|E_T|=m>1$. Let $L \in V_T$ be a leaf of the tree $T$, and let $P \in V_T$ be its unique parent. If we denote by $T^*$ the subtree of $T$ obtained by removing the edge $(P,L)$, then
\begin{equation*}
    \mathbb{E}\left[ \prod_{(i,j) \in E_T} f_{ij}(||X_i-X_j||) \right] =     \mathbb{E}\left[ f_{PL}(||X_P-X_L||) \prod_{(i,j) \in E_{T^*}} f_{ij}(||X_i-X_j||) \right].
\end{equation*}
Conditioning on $X_P$ results in
\begin{align*}
    \mathbb{E}\left[ \prod_{(i,j) \in E_T} f_{ij}(||X_i-X_j||) \right] &= \mathbb{E}\left[\mathbb{E}\left[ f_{PL}(||X_P-X_L||) \prod_{(i,j) \in E_{T^*}} f_{ij}(||X_i-X_j||) \; \middle| \; X_P \right]\right] \\
    &= \mathbb{E}\left[f_{PL}(||X_L||)\right] \mathbb{E}\left[\mathbb{E}\left[\prod_{(i,j) \in E_{T^*}} f_{ij}(||X_i-X_j||) \; \middle| \; X_P \right]\right] \\
    &= \mathbb{E}\left[f_{PL}(||X_P - X_L||)\right] \mathbb{E}\left[\prod_{(i,j) \in E_{T^*}} f_{ij}(||X_i-X_j||) \right].
\end{align*}
Applying the induction hypothesis on $T^*$ concludes the proof.
\end{proof}

\begin{proof}[Proof of Theorem \ref{thm:subgraphtree}]
The number of copies of $H$ inside the $\GIRGn$ as general subgraph is 
\begin{equation}
\begin{split}
    \Nsub(H) &= \sum_{\bv \in V_k} \mathbbm{1}_{\{H \subseteq \GIRGn|_{\bv}\}}. \label{eq:N(H)bernoulli}
\end{split}
\end{equation}
The connection probabilities in $\GIRGn$ depend on the vertex weights $\bw = \{w_i\}_{i \in [k]}$ and positions $\bx = \{x_i\}_{i \in [k]}$. Therefore,
\begin{equation*}
    \Prob{H \subseteq \GIRGn|_{\bv}} = \mathbb{E}_{\bw,\bx} \left[ \prod_{(i,j) \in \mathcal{E}} p(w_{v_i},w_{v_j},x_{v_i},x_{v_j}) \right].
\end{equation*}
Observe that $p(w_{v_i},w_{v_j},x_{v_i},x_{v_j})$ is a function of the distance $||x_{v_i} - x_{v_j}||$. Thus, for each $(i,j) \in E$ we can define $f_{ij}(||x_{v_i}-x_{v_j}||) := p(w_{v_i},w_{v_j},x_{v_i},x_{v_j})$. Then, the marginal probability given $\bw$ that $H\subseteq \GIRGn|_{\bv}$ is 
\begin{align*}
    \Prob{H\subseteq \GIRGn|_{\bv}} &= \mathbb{E}_{\bx} \left[ \prod_{(i,j) \in \mathcal{E}} f_{ij}(||x_{v_i}-x_{v_j}||) \right] \nonumber\\
    &= \prod_{(i,j) \in \mathcal{E}} \mathbb{E}_{\bx} \left[ f_{ij}(||x_{v_i}-x_{v_j}||) \right],
\end{align*}
where the last equality follows from Lemma \ref{lemma:treeedgesindependence}, as $H$ is a tree.
Moreover, by \cite[Lemma 3.3]{bringmann2019},
\begin{equation}\label{eq:marginalprobGIRG}
    \tilde{p}_{ij} := \mathbb{E}_{\bx} \left[ f_{ij}(||x_{v_i}-x_{v_j}||) \right] = \Theta\left(\frac{w_{v_i}w_{v_j}}{\mu n} \wedge 1\right)
\end{equation}
and combining \eqref{eq:marginalprobGIRG} with \eqref{eq:N(H)bernoulli},
\begin{equation*}
    \Nsub(H) \stackrel{d}{\sim}\sum_{\bv \in V_k} \prod_{(i,j) \in \mathcal{E}} \textnormal{Be}\left(\tilde{p}_{ij} \right).
\end{equation*}
This concludes the proof, as the right hand side of the latter equation is a random variable counting the occurrences of $H$ as general subgraph inside the IRG with weights $\bw$ and connection probabilities $\tilde{p}_{ij}$.
\end{proof}

\subsection{Proof for Hamiltonian patterns}
Finally, we prove Theorem \ref{thm:hamiltonian}. We will study the geometric and non-geometric regimes $k<1+\frac{1}{3-\tau}$ and $k>\frac{2}{3-\tau}$ separately. The proof strategy is similar in the two cases. First we show that the expectations of $\Nsub(H)$ and $\Nind(H)$ converge, when properly renormalized. The main challenge is to show that when a subgraph is Hamiltonian, this renormalized expectation is finite. Next, we will show that the two random variables concentrate around their means, which will lead to the main result.

\subsubsection*{Geometric regime $k<1+\frac{1}{3-\tau}$}

For any $H=(\mathcal{V},\mathcal{E})$, the expectation of $\Nsub(H)$ is
\begin{align}
    \mathbb{E}\left[\Nsub(H)\right] 
    &=  \frac{n!}{(n-k)!} \int \mathbb{P}(H \subseteq \GIRG(\bw,\bx)|_{\bv}) \dd \mathbb{P}(\bw,\bx) \nonumber\\
    &= \frac{n!}{(n-k)!} \int_{[w_0,\infty)^k} \prod_{i=1}^k f_W(w_{v_i}) \label{eq:integralform} \\
    &\hspace{2cm}\times \int_{(\mathbb{T}^d)^k} \prod_{(i,j) \in \mathcal{E}} \left(\frac{w_{v_i} w_{v_j}}{n\mu||x_{v_i}-x_{v_j}||^{d}} \wedge 1 \right)^{\gamma} \dd \bx  \dd \bw, \nonumber 
\end{align}
where $f_W$ is the probability distribution function of the power law random variable $W$ in~\eqref{eq:weightdistribution}.
The spatial coordinates $x_{v_1},...,x_{v_k}$ take values in $\mathbb{T}^d$. However, we may assume that $v_1$ has a fixed position $x_{v_1}$ (w.l.o.g. the origin of the torus) and obtain
\begin{align}
    \mathbb{E}\left[\Nsub(H)\right] &= \frac{n!}{(n-k)!} \int_{[w_0,\infty)^k} \prod_{i=1}^k f_W(w_{v_i}) \nonumber \\ & \hspace{1.5cm} \times \left(\int_{(\mathbb{T}^d)^{(k-1)}} \prod_{(i,j) \in \mathcal{E}} \left(\frac{w_{v_i} w_{v_j}}{n\mu||x_{v_i}-x_{v_j}||^{d}} \wedge 1 \right)^{\gamma} \dd x_{v_2} \cdots \dd x_{v_k} \right) \dd \bw \nonumber \\
    &= \frac{n!}{(n-k)!}\frac{1}{n^{k-1}} \int_{[w_0,\infty)^k} \prod_{i=1}^k f_W(w_{v_i})  \nonumber \\ 
    & \hspace{1cm} \times\left(\int_{(n^{1/d}\mathbb{T})^{d(k-1)}} \prod_{(i,j) \in \mathcal{E}} \left(\frac{w_{v_i} w_{v_j}}{\mu||\tilde{x}_{v_i}-\tilde{x}_{v_j}||^{d}} \wedge 1 \right)^{\gamma} \dd \tilde{x}_{v_2} \cdots \dd \tilde{x}_{v_k} \right) \dd \bw \nonumber \\
    &=: \frac{n!}{(n-k)!}\frac{1}{n^{k-1}} \Jsub_n(H) \label{eq:integralNsub}
\end{align}
where we have used the substitution $\tilde{x}_{v_i} = x_{v_i}/n^{1/d}$ for $i = 2,...,k$, and $n^{1/d}\mathbb{T}$ is the torus rescaled by a factor $n^{1/d}$.

In particular, as $n$ goes to infinity, $\Jsub_n(H)$ converges to the integral
\begin{multline}\label{eq:ConstGeoGeneral}
    \Jsub(H):= \int_{[w_0,\infty)^k} \prod_{i=1}^k f_W(w_{v_i}) \int_{(\mathbb{R}^d)^{k-1}} \prod_{(i,j) \in \mathcal{E}} \left(\frac{w_{v_i} w_{v_j}}{\mu||z_{v_i}-z_{v_j}||_{\infty}^{d}} \wedge 1 \right)^{\gamma} \\ \dd z_{v_2} \cdots \dd z_{v_k} \dd \bw,
\end{multline}
whenever $\Jsub(H)$ is defined.

In a similar way, we derive the expected number of induced subgraphs,
\begin{align}
    \mathbb{E}\left[\Nind(H)\right] &= \frac{n!}{(n-k)!}\frac{1}{n^{k-1}} \Jind_n(H). \label{eq:integralNind}
\end{align}
where $\Jind_n(H)$ is an integral converging to
\begin{multline}\label{eq:ConstGeoInduced}
    \Jind(H):= \int_{[w_0,\infty)^k} \prod_{i=1}^k f_W(w_{v_i}) \int_{(\mathbb{R}^d)^{k-1}} \prod_{(i,j) \in \mathcal{E}} \left(\frac{w_{v_i} w_{v_j}}{\mu||z_{v_i}-z_{v_j}||_{\infty}^{d}} \wedge 1 \right)^{\gamma} \\
    \times \prod_{(i,j) \not \in \mathcal{E}} 1 - \left(\frac{w_{v_i} w_{v_j}}{\mu||z_{v_i}-z_{v_j}||_{\infty}^{d}} \wedge 1 \right)^{\gamma} \dd z_{v_2} \cdots \dd z_{v_k} \dd \bw, 
\end{multline}
provided that $\Jind(H)$ is defined. We need to prove that when $k < 1 + 1/(3-\tau)$, 
\begin{equation*}
    \frac{\Nsub(H)}{n} \stackrel{p}{\longrightarrow} \Jsub(H),\qquad
    \frac{\Nind(H)}{n} \stackrel{p}{\longrightarrow} \Jind(H).
\end{equation*}
Then, we have all tools to prove convergence for the expectations of $\Nsub(H)/n$, $\Nind(H)/n$. 
\begin{proposition}\label{prop:ConvergenceMeanGeo}
    Under the hypothesis of Theorem \ref{thm:hamiltonian}(i), we have
    \begin{equation*}
        \frac{\E{\Nsub(H)}}{n} \longrightarrow \Jsub(H), \qquad \frac{\E{\Nind(H)}}{n} \longrightarrow \Jind(H)
    \end{equation*}
    where $\Jsub(H)$ and $\Jind(H)$ are the integrals defined in \eqref{eq:ConstGeoGeneral} and \eqref{eq:ConstGeoInduced}.
\end{proposition}

\begin{proof}
From equation \eqref{eq:expNind}, we can write
\begin{align}
    \mathbb{E}[\Nsub(H)] &= \frac{n!}{(n-k)!} \int \mathbb{P}(H \subset \GIRG(\bw,\bx)|_{\bv}) \dd \mathbb{P}(\bw,\bx) \nonumber \\
    &= \frac{n!}{(n-k)!} \int_{[w_0,\infty)^k} \prod_{i=1}^k f_W(w_{v_i}) \int_{(\mathbb{T}^d)k} \prod_{(i,j) \in \mathcal{E}} p(w_{v_i},w_{v_j},x_i,x_j) \dd \bx \dd \bw \nonumber 
\end{align}
If $H$ is a graph of size $k$ that contains a Hamiltonian cycle, we can assume (up to a permutation of the vertices) that $H$ contains the Hamiltonian cycle $\Hk$ with edge set $\mathcal{E}_{\Hk}=\{(i,i+1)\}_{i \in [k]}$. Moreover, introducing a factor $k$, we may assume that $w_{v_1}$ is the largest weight in $\bw$. Then, we upper bound the probability for the Hamiltonian pattern $H$ to appear on the vertex set $\bv$ with the probability that the open path from node $v_1$ to node $v_k$ appears:
\begin{align*}
    \mathbb{E}[\Nsub(H)] &= \frac{n!}{(n-k)!} \int_{[w_0,\infty)^k} \prod_{i=1}^k f_W(w_{v_i}) \int_{(\mathbb{T}^d)k} \prod_{(i,j) \in \mathcal{E}} p(w_{v_i},w_{v_j},x_i,x_j) \dd \bx \dd \bw \nonumber \\
    &\leq \frac{n!k}{(n-k)!} \int_{[w_0,\infty) \times [w_0,w_{v_1})^{k-1}} \prod_{i=1}^k f_W(w_{v_i}) \\& \hspace{4cm} \times \int_{(\mathbb{T}^d)k} \prod_{i \in [k-1]} p(w_{v_i},w_{v_{i+1}},x_i,x_{i+1}) \dd \bx \dd \bw \nonumber
\end{align*}
Now, from Lemma \ref{lemma:treeedgesindependence} we can average out the position from the connection probability. Since $\mathbb{E}_{\bx} \left[ p_{ij} \right] = \Theta\left(\frac{w_{v_i}w_{v_j}}{\mu n} \wedge 1\right)$, we have
\begin{equation}\label{eq:upperboundHk}
    \begin{aligned}
    \mathbb{E}[\Nsub(H)] &= O\left( \frac{n!k}{(n-k)!} \int_{[w_0,\infty) \times [w_0,w_{v_1})^{k-1}} \prod_{i=1}^k f_W(w_{v_i}) \prod_{i \in [k-1]} \frac{w_{v_i}w_{v_{i+1}}}{\mu n} \dd \bw \right) \\
    &= O \left( n^k \int_{[w_0,\infty) \times [w_0,w_{v_1})^{k-1}} \prod_{i=1}^k w_{v_i}^{-\tau} \prod_{i \in [k-1]} \frac{w_{v_i}w_{v_{i+1}}}{\mu n} \dd \bw \right)  \\
    &= O \left( n \int_{[w_0,\infty) \times [w_0,w_{v_1})^{k-1}} w_{v_1}^{1-\tau} w_{v_2}^{2-\tau} \cdots w_{v_{k-1}}^{2-\tau} w_{v_k}^{1-\tau} \dd \bw \right)  \\
    &= O \left( n \int_{[w_0,\infty)} w_{v_1}^{1-\tau + (3-\tau)(k-2)} \dd w_{v_1} \right). 
    \end{aligned}
\end{equation}
The latter integral converges whenever the exponent is smaller than $-1$, which is the case since $k < 1 + \frac{1}{3-\tau}$ by hypothesis.

On the other hand, $\mathbb{E}[\Nsub(H)] \geq \mathbb{E}[\Nsub(H;\Mab)]$, for any $\varepsilon, \balpha, \bbeta$. Setting $\balpha \equiv 0$ and $\bbeta \equiv - 1/d$, we have from \eqref{eq:Mabexponent}
\begin{equation}\label{eq:lowerboundHk}
    \mathbb{E}[\Nsub(H)] = \Omega(n)
\end{equation}
Now, combining \eqref{eq:upperboundHk}--\eqref{eq:lowerboundHk} yields $\mathbb{E}[\Nsub(H)] = \Theta(n)$, that is,

\begin{equation*}
    \limsup_{n \to \infty} \frac{\mathbb{E}[\Nsub(H)]}{n} < \infty, \qquad \liminf_{n \to \infty} \frac{\mathbb{E}[\Nsub(H)]}{n} > 0.
\end{equation*}
Then, from equation \eqref{eq:integralNsub}, 
\begin{align*}
    \lim_{n \to \infty} \frac{\mathbb{E}[\Nsub(H)]}{n} &= \lim_{n \to \infty} \frac{n!}{n^k(n-k)!} \Jsub_n(H) \nonumber \\ &= \lim_{n \to \infty} \Jsub_n(H) = \Jsub(H)\in (0, \infty). \nonumber
\end{align*}

Finally, in order to prove convergence for the expected number of induced Hamiltonian patterns, we first observe that $\mathbb{E}[\Nind(H)] \leq \mathbb{E}[\Nsub(H)]$. Indeed, if $H = \GIRG|_{\bv}$ for some $\bv$, then also $H \subseteq \GIRG|_{\bv}$. Moreover, since $(\balpha,\bbeta) \equiv (0,-1/d)$ is in the feasible region for the optimization problem \eqref{eq:optproblem}--\eqref{eq:optproblem_extra_constraint}, the expected number of copies of $H$ as induced subgraph is asymptotically lower bounded by $\Omega(n)$. Then, again, since $\mathbb{E}[\Nind(H)] = \Theta(n)$, we conclude that the integral $\Jind(H)$ is well defined, using the same argument as for $\mathbb{E}[\Nsub(H)]$ above, but replacing $\Jsub_n(H)$, $\Jsub(H)$ with $\Jind_n(H)$, $\Jind(H)$.\\
\end{proof}

\subsubsection*{Non-geometric regime $k > \frac{2}{3-\tau}$}
Using a similar approach as done in the geometric regime, we derive an explicit formula for $\E{\Nsub(H)}$ and $\E{\Nind(H)}$. From \eqref{eq:integralform}, using the substitution $y_{v_i} = w_{v_i}/\sqrt{n}$ for all $i = 1,...,k$, we have
\begin{align}
    \mathbb{E}\left[\Nsub(H)\right] 
    &= \frac{n!}{(n-k)!} \int_{[w_0,\infty)^k} \prod_{i=1}^k f_W(w_{v_i}) \nonumber \\
    &\hspace{2.5cm} \times \int_{(\mathbb{T}^d)^k} \prod_{(i,j) \in \mathcal{E}} \left(\frac{w_{v_i} w_{v_j}}{n\mu||x_{v_i}-x_{v_j}||^{d}} \wedge 1 \right)^{\gamma} \dd \bx \dd \bw \nonumber \\
    &= \frac{n!}{(n-k)!}n^{k(1-\tau)/2} \int_{[n^{-1/2}w_0,\infty)^k} \prod_{i=1}^k y_{v_i}^{-\tau} \nonumber \\
    &\hspace{2.5cm} \times \int_{(\mathbb{T}^d)^k} \prod_{(i,j) \in \mathcal{E}} \left(\frac{y_{v_i} y_{v_j}}{\mu||x_{v_i}-x_{v_j}||^{d}} \wedge 1 \right)^{\gamma} \dd \bx \dd \by \nonumber \\
    &=: \frac{n!}{(n-k)!}n^{k(1-\tau)/2} \Isub_n(H) \label{eq:integralNsubnongeo}
\end{align}
As $n \to \infty$, the integral $\Isub_n$ converges to 
\begin{equation}\label{eq:ConstNongeoGeneral}
    \Isub(H) := \int_{(0,\infty)^k} \prod_{i=1}^k y_{v_i}^{-\tau} \int_{(\mathbb{T}^d)^k} \prod_{(i,j) \in \mathcal{E}} \left(\frac{y_{v_i} y_{v_j}}{\mu||x_{v_i}-x_{v_j}||^{d}} \wedge 1 \right)^{\gamma} \dd \bx \dd \by,
\end{equation}
whenever $\Isub$ is defined. Similarly, we derive from \eqref{eq:expNind} that 
\begin{align*}
    \mathbb{E}\left[\Nind(H)\right] &= \frac{n!}{(n-k)!}n^{k(1-\tau)/2} \Iind_n(H),
\end{align*}
where $\Iind_n(H)$ is an integral converging to
\begin{equation}\label{eq:ConstNongeoInduced}
    \begin{aligned}
    \Iind(H) &:= \int_{(0,\infty)^k} \prod_{i=1}^k y_{v_i}^{-\tau} \int_{(\mathbb{T}^d)^k}  \prod_{(i,j) \in \mathcal{E}} \left(\frac{y_{v_i} y_{v_j}}{\mu||x_{v_i}-x_{v_j}||^{d}} \wedge 1 \right)^{\gamma} \\
    &\hspace{3.5cm} \times \prod_{(i,j) \not \in \mathcal{E}} 1 - \left(\frac{y_{v_i} y_{v_j}}{\mu||x_{v_i}-x_{v_j}||^{d}} \wedge 1 \right)^{\gamma} \dd \bx \dd \by, 
    \end{aligned}
\end{equation}
provided that $\Iind(H)$ is defined. We now prove convergence for the expectations of $\Nsub(H)/n^{k(3-\tau)/2}$ and $\Nind(H)/n^{k(3-\tau)/2}$.

\begin{proposition}\label{prop:ConvergenceMeanNongeo}
    Under the hypothesis of Theorem \ref{thm:hamiltonian}(ii), we have
    \begin{equation*}
        \frac{\E{\Nsub(H)}}{n^{k(3-\tau)/2}} \longrightarrow \Isub, \qquad \frac{\E{\Nind(H)}}{n^{k(3-\tau)/2}} \longrightarrow \Iind
    \end{equation*}
    where $\Isub(H)$ and $\Iind(H)$ are the integrals defined in \eqref{eq:ConstNongeoGeneral} and \eqref{eq:ConstNongeoInduced}.
\end{proposition}

\begin{proof}
We start by proving that $\Isub(H)$ and $\Iind(H)$ are well defined. By the assumptions of Theorem \ref{thm:hamiltonian}, we can assume without loss of generality that $\Hk \subseteq H$. Then,
\begin{equation*}
    \Iind(H) \leq \Isub(H) \leq \Isub(\Hk).
\end{equation*}
The first inequality is obtained by upper bounding the product over $(i,j) \notin \mathcal{E}$ in  \eqref{eq:ConstNongeoInduced} with 1; whereas in the second inequality we upper bound all the terms in the product over $(i,j) \in \mathcal{E}$, when $(i,j) \not \in \mathcal{E}_{\Hk}$ in \eqref{eq:ConstNongeoGeneral} with 1.
Therefore, $\Isub(H)$ and $\Iind(H)$ are well defined as soon as $\Isub(\Hk) < \infty$. Since $k$ is odd, we can apply the substitution of variables $u_{v_i} = y_{v_i}y_{v_{i+1}}$ and rewrite $\Isub(\Hk)$ as

\begin{align*}
    \Isub(\Hk) &= \int_{(0,\infty)^k} \prod_{i=1}^k y_{v_i}^{-\tau} \int_{(\mathbb{T}^d)^k} \prod_{i = 1}^k \left(\frac{y_{v_i} y_{v_{i+1}}}{\mu||x_{v_i}-x_{v_{i+1}}||^{d}} \wedge 1 \right)^{\gamma} \dd \bx \dd \by \nonumber\\
    &= \int_{(0,\infty)^k} \prod_{i=1}^k u_{v_i}^{-(\tau+1)/2} \int_{(\mathbb{T}^d)^k} \prod_{i = 1}^k \left(\frac{u_{v_i}}{\mu||x_{v_i}-x_{v_{i+1}}||^{d}} \wedge 1 \right)^{\gamma} \dd \bx \dd \bu.
\end{align*}
Next, due to the symmetry of the Hamiltonian cycle $\Hk$, we assume without loss of generality that $u_{v_k}$ is the maximum among all $u_{v_i}$, for $i = 1,...,k$. Then

\begin{align*}
    \Isub(\Hk) &= k \int_{(0,u_{v_k})^{k-1}\times (0,\infty)} \prod_{i=1}^k u_{v_i}^{-(\tau+1)/2} \int_{(\mathbb{T}^d)^k} \prod_{i = 1}^k \left(\frac{u_{v_i}}{\mu||x_{v_i}-x_{v_{i+1}}||^{d}} \wedge 1 \right)^{\gamma} \dd \bx \dd \bu \nonumber \\
    &\leq k \int_{(0,u_{v_k})^{k-1}\times (0,\infty)} \prod_{i=1}^k u_{v_i}^{-(\tau+1)/2} \int_{(\mathbb{T}^d)^k} \prod_{i = 1}^{k-1} \left(\frac{u_{v_i}}{\mu||x_{v_i}-x_{v_{i+1}}||^{d}} \wedge 1 \right)^{\gamma} \dd \bx \dd \bu ,
\end{align*}
where in the last inequality we have used that $\left(\frac{u_{v_k}}{\mu||x_{v_k}-x_{v_{1}}||^{d}} \wedge 1 \right) \leq 1$. Now, since we are left with an integral over the probability of an open path, we can use Lemma \ref{lemma:treeedgesindependence} and approximate each connection probability by
\begin{equation*}
    \int_{(\mathbb{T}^d)^k} \prod_{i = 1}^{k-1} \left(\frac{u_{v_i}}{\mu||x_{v_i}-x_{v_{i+1}}||^{d}} \wedge 1 \right)^{\gamma} \dd \bx = \Theta\left( \prod_{i = 1}^{k-1} \left( u_{v_i} \wedge 1 \right) \right).
\end{equation*}
Then,
\begin{equation*}
    \Isub(\Hk) \leq k \int_{(0,u_{v_k})^{k-1}\times (0,\infty)} \prod_{i=1}^k u_{v_i}^{-(\tau+1)/2} \prod_{i = 1}^{k-1} \Theta\left( u_{v_i} \wedge 1 \right) \dd \bu .
\end{equation*}
This integrand is separable in all variables, and we distinguish two cases.
\begin{itemize}
    \item When $u_{v_k} \leq 1$, we have $(u_{v_i} \wedge 1) = u_{v_i}$ for all $i = 1,...,k-1$. Then, for all $i = 1,...,k-1$,
    \begin{equation*}
        \int_0^{u_{v_k}} (u_{v_i} \wedge 1) u_{v_i}^{-(\tau+1)/2} \dd u_{v_i} = \int_0^{u_{v_k}} u^{-(\tau-1)/2} \dd u = c_1 \cdot u_{v_k}^{(3-\tau)/2}. 
    \end{equation*}
    for some constant $c_1 > 0$.
    \item When $u_{v_k} > 1$, we split the integral over $u_{v_i}$ for each $i = 1,...,k-1$, into two separate regions
    \begin{align*}
        \int_0^{u_{v_k}} (u_{v_i} \wedge 1) u_{v_i}^{-(\tau+1)/2} \dd u_{v_i} &= \int_0^1 u^{-(\tau-1)/2} \dd u + \int_1^{u_{v_k}} u^{-(\tau+1)/2} \dd u \\
        &= c_2 + (c_3 - u_{v_k}^{-(\tau-1)/2}) \leq c_4,
    \end{align*}
    as $\tau\in(2,3)$, for some constants $c_2,c_3,c_4>0$.
\end{itemize}
Thus, integrating $u_{v_k}$ over the two different regions $(0,1)$ and $(1,\infty)$ gives:
\begin{align*}
    \Isub(\Hk) &\leq k \int_{0}^{\infty} \left( \iint_0^{u_{v_k}} \prod_{i=1}^k u_{v_i}^{-(\tau+1)/2} \prod_{i = 1}^{k-1} \left( u_{v_i} \wedge 1 \right) \dd u_{v_1} \cdots \dd u_{v_{k-1}} \right) \dd u_{v_k} \\
    &= k \left( \int_{0}^{1} c_1 \cdot u^{\frac{3-\tau}{2}(k-1) - \frac{\tau+1}{2}} \dd u + \int_{1}^{\infty} c_4 \cdot u^{- \frac{\tau+1}{2}} \dd u \right).
\end{align*}
The integral over $[1,\infty)$ is finite, since $-\frac{\tau + 1}{2} < - 1$. Moreover, the integral over $[0,1]$ is finite as soon as $\frac{3-\tau}{2}(k-1) - \frac{\tau+1}{2} > - 1$, which corresponds exactly to the hypothesis $k > \frac{2}{3-\tau}$ of Theorem \ref{thm:hamiltonian}(ii). Hence, since $\Isub(\Hk) < \infty$, both integrals $\Isub(H)$ and $\Iind(H)$ are well defined.

Finally, from \eqref{eq:integralNsubnongeo} we conclude that
\begin{align*}
    \lim_{n \to \infty} \frac{\E{\Nsub(H)}}{n^{k(3-\tau)/2}} &= \lim_{n \to \infty} \frac{\frac{n!}{(n-k)!}n^{k(1-\tau)/2} \Isub_n(H) }{ n^{k(3-\tau)/2}} \\
    &= \lim_{n \to \infty} \Isub_n(H) = \Isub(H) < \infty,
\end{align*}
and similarly for $\Nind(H)$.
\end{proof}

\subsubsection*{Concentration around the mean}
We have now all tools to demonstrate the concentration result for Hamiltonian patterns. 
\begin{proof}[Proof of Theorem \ref{thm:hamiltonian}(i)]
In this case, $k< 1 + 1/(3-\tau)$. Set $\star \in \{\text{sub},\text{ind}\}$ and $(\balpha^*,\bbeta^*)\equiv(0,-1/d)$. Observe that, for any $\varepsilon > 0$,
\begin{equation*}
    \mathbb{E}[\Nstar(H)] = \mathbb{E}[\Nstar(H;\M)] + \mathbb{E}[\Nstar(H;\Mbar)],
\end{equation*}
where $\M$ is defined as in \eqref{eq:Mab} and $\Mbar$ is the complementary set. Moreover, there exists $J_{\varepsilon}(H)$ such that
\begin{align*}
    \mathbb{E}[\Nstar(H;\M)] = n J_{\varepsilon}(H)(1+o(1)), \quad \text{ and } \quad J_{\varepsilon}(H) \to \Jstar(H) \quad \text{as } \varepsilon \to 0.
\end{align*}
by Proposition \ref{prop:ConvergenceMeanGeo} combined with the fact that the sequence of random variables $\Nstar(H,\M)$ indexed by $\varepsilon$ converges almost surely to $\Nstar(H)$. Now, by Lemma \ref{lemma:selfaveraging} and the Chebychev inequality we can conclude convergence in probability in the optimal set:
\begin{equation}\label{eq:probconvMeps}
    \frac{\Nstar(H;\M)}{n} \stackrel{p}{\longrightarrow} J_{\varepsilon}(H)
\end{equation}
On the other hand, the Markov inequality shows that $\Nstar(H;\Mbar)$ is asymptotically upper bounded (in probability) by $(\Jstar(H)- J_{\varepsilon}(H)) \cdot O(n)$. Then, for any $\delta>0$,
\begin{equation}\label{eq:probboundMeps}
    \limsup_{\varepsilon \to 0} \limsup_{n \to \infty} \mathbb{P}\left(\frac{\Nstar(H;\Mbar)}{n} > \delta \right) = 0
\end{equation}
Finally, combining \eqref{eq:probconvMeps} and \eqref{eq:probboundMeps} we conclude that
\begin{equation*}
    \frac{\Nstar(H)}{n} \stackrel{p}{\longrightarrow} \Jstar(H).
\end{equation*}
\end{proof}

\begin{proof}[Proof of Theorem \ref{thm:hamiltonian}(ii)]
The proof follows the same steps as in (i). In this case, since $k$ is odd and $k>2/(3-\tau)$, Proposition \ref{prop:ConvergenceMeanNongeo} holds, and to obtain expressions similar to  \eqref{eq:probconvMeps} and \eqref{eq:probboundMeps} we need to renormalize by $n^{k(3-\tau)/2}$ and replace $\Jstar(H)$ with $\Istar(H)$.
\end{proof}

%
%

\begin{funding}
The second author was supported by NWO VENI grant 202.001
\end{funding}



\bibliographystyle{imsart-number} 
\bibliography{sample}       


\end{document}